\documentclass[10pt]{article}
\pdfoutput=1
\usepackage[square,authoryear]{natbib}
\usepackage{graphicx}
\usepackage{wrapfig}
\usepackage{tikz}
\usepackage{url}
\usepackage{marsden_article}
\newtheorem{exams}[theorem]{Examples}
\newtheorem{remark}[theorem]{Remark}
\newtheorem{assumption}[theorem]{Assumption}

\newcommand{\cC}{\mathcal{C}}

\begin{document}

\title{Discrete Mechanics and Optimal Control: an Analysis}

\author{Sina Ober-Bl\"obaum\thanks {Research partially supported by 
the University of Paderborn, Germany and AFOSR grant FA9550-08-1-0173; sinaob@cds.caltech.edu}
 \\ Control and Dynamical Systems
\\ California Institute of Technology 107-81
\\ Pasadena, CA 91125\\
\and
Oliver Junge\thanks{email: junge@ma.tum.de}\\
Zentrum Mathematik\\
Technische Universit\"{a}t MŸ\"{u}nchen\\
D-85747 Garching, Germany \\
 \and 
 Jerrold E. Marsden\thanks{Research partially supported by 
AFOSR grant FA9550-08-1-0173; email: jmarsden@caltech.edu}
\\ Control and Dynamical Systems
\\ California Institute of Technology 107-81
\\ Pasadena, CA 91125
}

\date{\small Version 0.4: Sep 8, 2008}

\maketitle

\begin{abstract}
The optimal control of a mechanical system is of crucial importance in many realms. Typical examples are the determination of a time-minimal path in vehicle dynamics, a minimal energy trajectory in space mission design, or optimal motion sequences in robotics and biomechanics.  In most cases, some sort of discretization of the original, infinite-dimensional optimization problem has to be performed in order to make the problem amenable to computations.  The approach proposed in this paper is to directly discretize the variational description of the system's motion. The resulting optimization algorithm lets the discrete solution directly inherit characteristic structural properties from the continuous one like symmetries and integrals of the motion. We show that the DMOC approach is equivalent to a finite difference discretization of Hamilton's equations by a symplectic partitioned Runge-Kutta scheme and employ this fact in order to give a proof of convergence. 
 The numerical performance of DMOC and its relationship to other existing optimal control methods are investigated.  
\end{abstract}

\tableofcontents

\section{Introduction}

In order to solve optimal control problems for mechanical systems, this paper links two important areas of research: {\em optimal control} and {\em variational mechanics}. The motivation for combining these fields of investigation is twofold. Besides the aim of preserving certain properties of the mechanical system for the approximated optimal solution, optimal control theory and variational mechanics have their common origin in the {\em calculus of variations}. In mechanics, the calculus of variations is also fundamental through the {\em principle of stationary action} also called {\it Hamilton's principle}. When applied to the action of a mechanical system, this principle yields the equations of motion for that system---the {\em Euler-Lagrange equations}.  In optimal control theory the  calculus of variations also plays a fundamental role. For example, it is used to derive optimality conditions via the {\em Pontryagin maximum principle}. In addition to its importance in continuous mechanics and control theory, the {\em discrete calculus of variations} and the corresponding discrete variational principles play an important role in constructing efficient numerical schemes for the simulation of mechanical systems and for optimizing dynamical systems.

\paragraph{Discrete Optimal Control and Discrete Variational Mechanics.}
The theory of discrete variational mechanics has its roots in the optimal control literature of the $1960$s; see for example \cite{Cad70, HwFa67, JoPo64}. Specifically, \cite{Cad70} developed a discrete calculus of variations theory in the following way: A function is introduced which depends on a sequence of numbers, e.g.\ a sequence of times. A minimizing sequence necessarily satisfies a second-order difference equation, which is called the {\em discrete Euler equation} in reminiscence of its similarity with the Euler equation of the classical calculus of variations.

An application of the discrete calculus of variations to an optimal control problem leads to a so called {\em direct solution method}. In this, one transforms the optimal control problem into a finite dimensional equality constrained nonlinear optimization problem via a finite dimensional parametrization of states and controls.   In contrast, \emph{indirect methods} (see \S\ref{corres_ocp} for an overview) are based on the explicit derivation of the necessary optimality conditions via the Pontryagin maximum principle. 

On the other hand, the theory of discrete variational mechanics describes a variational approach to discrete mechanics and mechanical integrators. The application of a discrete version of Hamilton's principle results in the {\em discrete Euler-Lagrange equations}. Analogous to the continuous case, near conservation of discrete energy, discrete momentum maps related to the discrete system's symmetries and the discrete symplectic form can be shown. This is due to the discretization of the variational structure of the mechanical system directly. Early work on discrete mechanics was often independently done by \cite{Cad73, Lee83, Lee87, Log73, Maeda80, Maeda81a, Maeda81b}. In this work, the role of the discrete action sum, the discrete Euler-Lagrange equations and the discrete Noether's theorem were clearly understood. The variational view of discrete mechanics and its numerical implementation is further developed in \cite{WM97,WM97b} and then extended in \cite{BoSu99a,BoSu99b,Kane99,Kane00,MaPeSh99a,MaPeSh99b}. The route of a variational approach to symplectic-momentum integrators has been taken by \cite{Sur90}, \cite{Mac92}; see the review by  \cite{MaWe01} and references therein. In this review a detailed derivation and investigation of these {\em variational integrators} for conservative as well as for forced and constrained systems is given. 

\paragraph{Combining Optimal Control and Variational Mechanics.} The present paper concerns the optimal control of dynamical systems whose behavior can be described by the {\em Lagrange-d'Alembert principle}.  To numerically solve this kind of problem, we make use of the discrete calculus of variations only, that means we apply the discrete variational principle on two layers. On the one hand we use it for the description of the mechanical system under consideration, and on the other hand for the derivation of necessary optimality conditions for the optimal control problem.  The application of discrete variational principles already on the dynamical level (namely the discretization of the Lagrange-d'Alembert principle) leads  to structure-preserving time-stepping equations which serve as equality constraints for the resulting finite dimensional nonlinear optimization problem. The benefits of variational integrators are handed down to the optimal control context. For example, in the presence of symmetry groups in the continuous dynamical system, also along the discrete trajectory the change in momentum maps is consistent with the control forces. Choosing the objective function to represent the control effort, which has to be minimized is only meaningful if the system responds exactly according to the control forces. 

\paragraph{Related Work.} A survey of different methods for the optimal control of dynamical systems described by ordinary differential equations is given in \S\ref{corres_ocp}. However, to our knowledge, DMOC is the first approach to solutions of optimal control problems involving the concept of discrete mechanics to derive structure-preserving schemes for the resulting optimization algorithm. Since our first formulations and applications to space mission design and formation flying (\cite{JMOb,JMaOb,JOb,OB08}), DMOC has been applied for example to problems from robotics and biomechanics (\cite{KaMa05, KDMS07, KoSu07, JM06, pekam07, ross05}) and to image analysis (\cite{McLachlan06a}). From the theoretical point of view, considering the development of variational integrators, extensions of DMOC to mechanical systems with nonholonomic constraints or to systems with symmetries are quite natural and have already been analyzed in \cite{KDMS07,KoSu07}. Further extensions are currently under investigation, for example DMOC for hybrid systems \cite{pekam07} and for constrained multi-body dynamics (see \cite{DMOCC,LObMO07,OB08}). DMOC related approaches are presented in \cite{LML2, LML1}. The authors discretize the dynamics by a Lie group variational integrator. Rather than solving the resulting optimization problem numerically, they construct the discrete necessary optimality conditions via the discrete variational principle and solve the resulting discrete boundary value problem (the discrete state and adjoint system). The method is applied to the optimal control of a rigid body and to the computation of attitude maneuvers of a rigid spacecraft. 

\paragraph{Outline.} In \S\ref{varmech} and \S\ref{discrete_forcing} we introduce the relevant concepts from classical variational mechanics and discrete variational mechanics following the work of \cite{MaWe01}. Especially, we focus on the Lagrangian and Hamiltonian description of control forces for the established framework of variational mechanics. Definitions and concepts of the variational principle, the Legendre transform, and Noether's theorem are readopted for the forced case in both the continuous as well as the discrete setting. 
In \S\ref{sec:locp} and \S\ref{sec:oc_dm} we combine concepts from optimal control and discrete variational mechanics to build up a setting for the optimal control of a continuous and a discrete mechanical system, respectively.
\S\ref{corres_cd} describes the correspondence between the continuous and the discrete Lagrangian system as basis for a comparison between the continuous and the discrete optimal control problems in \S\ref{corres_ocp}:  We link both frameworks viewing the discrete problem as an approximation of the continuous one. 
The application of discrete variational principles for a discrete description of the dynamical system leads to structure-preserving time-stepping 
equations. Here, the special benefits of variational integrators are handed down to the optimal control context. 
These time-stepping equations serve as equality constraints for the resulting finite dimensional nonlinear optimization problem, therefore the described procedure can be categorized as a direct solution method.
Furthermore, we show the equivalence of the discrete Lagrangian optimal control problems to those resulting from Runge-Kutta discretizations of the corresponding Hamiltonian system. This equivalence allows us to construct and compare the adjoint systems of the continuous and the discrete Lagrangian optimal control problem. In this way, one of our main results is related to the order of approximation of the adjoint system of the discrete optimal control problem to that of the continuous one. 
With the help of this approximation result, we show that the solution of the discrete Lagrangian optimal control system converges to the continuous solution of the original optimal control problem. The proof strategy is based on existing convergence results of optimal control problems discretized via Runge-Kutta methods (\cite{DoHa00,Hager00}).
\S\ref{impl} gives a detailed description of implementation issues of our method. Furthermore, in \S\ref{appl} we numerically verify the preservation and convergence properties of DMOC and the benefits of using DMOC compared to other standard methods to the solution of optimal control problems. 

\section{Mechanical systems with forcing and control}

\subsection{Variational Mechanics}\label{varmech}

Our aim is to optimally control Lagrangian and Hamiltonian systems.  For the description of their dynamics, we introduce a variational framework including external forcing resulting from dissipation, friction, loading and in particular control forces. To this end, we extend the notions in \cite{MaWe01} to Lagrangian control forces.

\paragraph{Forced Lagrangian Systems.} Consider an $n$-dimensional configuration manifold $Q$ with local coordinates $q=(q^1,\dots,q^n)$, the associated state space given by the tangent bundle $TQ$ and a Lagrangian $L: TQ \rightarrow \mathbb{R}$. 
Given a time interval $[0,T]$, we define the \textit{path space}
\[ \mathcal{C} (Q) = \mathcal{C}([0,T],Q) = \{q:[0,T] \rightarrow Q \, | \, q\;\text{is a $C^2$ curve}\} \]
and the \textit{action map} $\mathfrak{G}: \mathcal{C} \rightarrow \mathbb{R}$
\begin{equation}\label{eq:actionmap}
\mathfrak{G}(q) = \int _0^T L(q(t),\dot{q}(t)) \, \text{d} t.
\end{equation}
The tangent space $T_q\mathcal{C}(Q)$ to $\mathcal{C}(Q)$ at the point $q$ is the set of $C^2$ maps $v_q: [0,T]\rightarrow TQ$ such that $\tau_Q\circ v_q = q$, where $\tau_Q:TQ \rightarrow Q$ is the natural projection.


To define control forces for Lagrangian systems, we introduce a \textit{control manifold} $U\subset \mathbb{R}^m$ and define the \textit{control path space} 
\[
\mathcal{C}(U) = \mathcal{C}([0,T],U) = \{ u: [0,T] \rightarrow U\,|\, u \in L^\infty\},
\]
with $u(t) \in U$ also called the \textit{control parameter}. 
Here, $L^\infty$ denotes the space of essentially bounded, measurable functions equipped with the essential supremum norm.
With this notation we define a \index{Lagrangian control force|textit}{\em Lagrangian control force} as a map $f_{L}: TQ \times U \rightarrow T^*Q$, which is given in coordinates as
\[f_{L}: (q,\dot{q},u) \mapsto (q,f_{L}(q,\dot{q},u)),\]
where we assume that the control forces can also include configuration and velocity dependent forces resulting e.g.~from dissipation and friction.
We interpret a Lagrangian control force as a family of Lagrangian forces with fixed curves $u$, that are fiber-preserving maps $f_{L}^u: TQ  \rightarrow T^*Q$ over the identity $id_Q$, which we write in coordinates as
\[f_{L}^u: (q,\dot{q}) \mapsto (q,f_{L}^u(q,\dot{q})).\]
Whenever we denote $f_{L}(q,\dot{q},u)$ as a one-form on $TQ$, we mean the family of horizontal one-forms $f_{L}^u(q,\dot{q})$ on $TQ$ induced by the family of fiber-preserving maps $f_{L}^u$.

 
Given a control path $u\in\mathcal{C}(U)$, the \textit{Lagrange-d'Alembert principle} seeks curves $q\in \mathcal{C}(Q)$ satisfying
\begin{equation}\label{eq:varmech_ldap}
\delta \int _0^T L(q(t),\dot{q}(t)) \, \text{d} t + \int _0^T f_{L}(q(t),\dot{q}(t),u(t)) \cdot \delta q(t)\, \text{d} t = 0,
\end{equation}
where $\delta$ represents variations vanishing at the endpoints. The second integral in \eqref{eq:varmech_ldap} is the \textit{virtual work} acting on the mechanical system via the force $f_{L}$. Integration by parts 
shows that this is equivalent to the \textit{forced Euler-Lagrange equations}
\begin{equation}\label{eq:varmech_FEL}
\frac{\partial L}{\partial q}(q,\dot{q}) - \frac{\text{d}}{\text{d} t} \left( \frac{\partial L}{\partial \dot{q}}(q,\dot{q}) \right) + f_{L}(q,\dot{q},u) = 0.
\end{equation}
These equations implicitly define a family of \textit{forced Lagrangian flows} and \textit{forced Lagrangian vector fields} where for a fixed curve $u\in\mathcal{C}(U)$, $F_L^u: TQ \times \mathbb{R} \rightarrow TQ$ denotes the forced Lagrangian flow of the forced Lagrangian vector field $X_L^u: TQ \rightarrow T(TQ)$.

The one-form $\Theta_L$ on $TQ$ given in coordinates  by
\begin{equation}\label{eq:oneform}
\Theta_L = \frac{\partial L}{\partial \dot{q}^i} \mathbf{d}  q^i
\end{equation}
is called the \textit{Lagrangian one-form} and
the \textit{Lagrangian symplectic form} $\Omega_L = \mathbf{d}  \Theta_L$,  is given in coordinates by
\[\Omega_L(q,\dot{q}) = \frac{\partial^2L}{\partial q^i \partial \dot{q}^j} \mathbf{d}  q^i \wedge \mathbf{d}  q^j + \frac{\partial^2L}{\partial \dot{q}^i \partial \dot{q}^j} \mathbf{d}  \dot{q}^i \wedge \mathbf{d}  q^j.\] 

Recall that the absence of forces, the Lagrangian symplectic form is preserved under the Lagrangian flow as $(F_L^T)^*(\Omega_L)= \Omega_L$, where $F_L^t:TQ\rightarrow TQ$ is the map $F_L$ at some fixed time $t$.

\paragraph{Forced Hamiltonian Systems.}

Consider an $n$-dimensional configuration manifold $Q$, and define the phase space to be the cotangent bundle $T^*Q$. The Hamiltonian is a function $H: T^*Q \rightarrow \mathbb{R}$. We will take local coordinates on $T^*Q$ to be $(q,p)$ with $q=(q^1,\dots,q^n)$ and $p=(p_1,\dots,p_n)$.
Define the \textit{canonical one-form} $\Theta$ on $T^*Q$ by
\begin{equation}\label{eq:varmech_canon1form}
\Theta(p_q)\cdot u_{p_q} = \langle p_q, T_{\pi_{Q}} \cdot u_{p_q} \rangle,
\end{equation}
where $p_q \in T^*Q$, $u_{p_q}\in T_{p_q}(T^*Q)$,  $\pi_{Q}: T^*Q \rightarrow Q$ is the canonical projection, $T_{\pi_{Q}}: T(T^*Q)\rightarrow TQ$ is the tangent map of $\pi_{Q}$ and $\langle \cdot, \cdot \rangle$ denotes the natural pairing between vectors and covectors. In coordinates, we have $\Theta(q,p)= p_i\mathbf{d}  q^i$. The \textit{canonical two-form} $\Omega$ on $T^*Q$ is defined to be
\[\Omega = -\mathbf{d}  \Theta,\]
which has the coordinate expression $\Omega(q,p) = \mathbf{d}  q^i \wedge  \mathbf{d}  p_i$.

A \textit{Hamiltonian control force} is a map $f_{H}: T^*Q \times U \rightarrow T^*Q$  identified by a family of fiber preserving maps $f_H^u: T^*Q \rightarrow T^*Q$ over the identity. Given such a control force, we define the corresponding family of horizontal one-forms $f'_H$ on $T^*Q$ by
\[(f^{u}_H)'(p_q) \cdot w_{p_q} = \langle f_H^u(p_q), T_{\pi_Q} \cdot w_{p_q}\rangle\]
for fixed curves $u\in\mathcal{C}(U)$ where $\pi_Q: T^*Q \rightarrow Q$ is the projection. This expression is reminiscent of definition \eqref{eq:varmech_canon1form} of the canonical one-form $\Theta$ on $T^*Q$, and in coordinates it reads $(f^u_H)'(q,p) \cdot (\delta q, \delta p) = f^u_{H}(q,p) \cdot \delta q$, thus the one-form is clearly horizontal.
The \textit{forced Hamiltonian vector field} $X^u_H$ for a fixed curve $u$ is now defined by the equation
\[\mathbf{i}  _{X^u_H} \Omega = \mathbf{d}  H - (f^u_H)'\]
and in coordinates this gives the well-known \textit{forced Hamilton's equations}
\begin{subequations}\label{eq:varmech_forcHam}
\begin{align}
X^u_q(q,p)& =  \frac{\partial H}{\partial p}(q,p),\label{eq:varmech_forcHam_q}\\
X^u_p(q,p) &=  -\frac{\partial H}{\partial q}(q,p) + f^u_{H}(q,p),\label{eq:varmech_forcHam_p}
\end{align}
\end{subequations}
which are the standard Hamilton's equations in coordinates with the forcing term added to the momentum equation. This defines the \index{forced Hamiltonian flow|textit}{\em forced Hamiltonian flow} $F^u_H: T^*Q \times \mathbb{R} \rightarrow T^*Q$ of the forced Hamiltonian vector field $X^u_H= (X^u_q, X^u_p)$ for a fixed curve $u\in \mathcal{C}(U)$.

\paragraph{The Legendre Transform with Forces.}
Given a Lagrangian $L$, we can take the standard \textit{Legendre transform} $\mathbb{F} L: TQ \rightarrow T^*Q$ defined by 
\[\mathbb{F} L (v_q) \cdot w_q = \left.\frac{\text{d}}{\text{d} \epsilon}\right|_{\epsilon=0} L(v_q + \epsilon w_q), \]
where $v_q, w_q \in T_qQ$, and which has coordinate form
\[\mathbb{F} L: (q,\dot{q}) \mapsto (q,p)= \left( q,\frac{\partial L}{\partial \dot{q}} (q,\dot{q})\right),\]
and relate Hamiltonian and Langrangian control forces by
\[ f_{L}^u = f_{H}^u \circ \mathbb{F} L.\]
If we also have a Hamiltonian $H$ related to $L$ by the Legendre transform as
\[ H(q,p) = \mathbb{F}L(q,\dot{q})\cdot \dot{q} - L(q,\dot{q}), \]
then the forced Euler-Lagrange equations and the forced Hamilton's equations are equivalent. That is, if $X_L^u$ and $X_H^u$ are the forced Lagrangian and Hamiltonian vector fields, respectively, then $(\mathbb{F} L)^*(X_H^u) = X_L^u$. To see this, we compute

\begin{align*}
\frac{\partial H}{\partial q}(q,p) & = p\cdot \frac{\partial \dot{q}}{\partial q}  - \frac{\partial L}{\partial {q}}(q,\dot{q}) - \frac{\partial L}{\partial \dot{q}} (q,\dot{q}) \frac{\partial \dot{q}}{\partial q} \\
& = - \frac{\partial L}{\partial {q}}(q,\dot{q})\\
& = -\frac{\text{d}}{\text{d} t} \left( \frac{\partial L}{\partial \dot{q}}(q,\dot{q})\right) + f_{L}^u(q,\dot{q})\\
& = -\dot{p} + f_{H}^u \circ \mathbb{F} L(q,\dot{q}),\\
& = -\dot{p} + f_{H}^u(q,p) ,\\
\frac{\partial H}{\partial p}(q,p) & = \dot{q} + p\cdot \frac{\partial \dot{q}}{\partial p} - \frac{\partial L}{\partial \dot{q}} (q,\dot{q}) \frac{\partial \dot{q}}{\partial p}\\
& = \dot{q}, 
\end{align*}
where $p = \mathbb{F} L (q,\dot{q})$  defines $\dot{q}$ as a function of $(q,p)$.

\paragraph{Noether's Theorem with Forcing.}
A key property of Lagrangian flows is their behavior with respect to group actions. Assume a Lie group $G$ with Lie algebra $\mathfrak{g}$ acts on $Q$ by the (left or right) action $\phi: G\times Q\rightarrow Q$. Consider the tangent lift of this action to $\phi^{TQ}:G\times TQ$ given by $\phi_g^{TQ}(v_q)=T(\phi_g)\cdot v_q$, which is
 \[\phi^{TQ}(g,(q,\dot{q})) = \left(\phi^i(g,q),\frac{\partial \phi^i}{\partial q^j}(g,q) \dot{q}^j  \right).\]
 For $\xi\in \mathfrak{g}$ the {\em infinitesimal generators} $\xi_{Q}:Q \rightarrow TQ$ and $\xi_{TQ}:TQ \rightarrow T(TQ)$ are defined by
\[\xi_{Q}(q) = \frac{\text{d}}{\text{d} g} \left( \phi_g(q) \right)\cdot \xi,\quad \xi_{TQ}(v_q) = \frac{\text{d}}{\text{d} g} \left( \phi_g^{TQ}(v_q) \right)\cdot \xi,\]
and the {\em Lagrangian momentum map} $\mathbf{J} _L:TQ\rightarrow \mathfrak{g}^*$  is defined to be
\[ \mathbf{J} _L(v_q)\cdot \xi = \Theta_L \cdot \xi_{TQ}(v_q).\]
If the Lagrangian is {\em invariant} under the lift of the action, that is we have $L \circ \phi_g^{TQ} = L$ for all $g\in G$ (we also say, the group action is a {\em symmetry} of the Lagrangian), the Lagrangian momentum map is preserved of the Lagrangian flow in the absence of external forces. We now consider the effect of forcing on the evolution of momentum maps that arise from symmetries of the Lagrangian.
In \cite{MaWe01} it is shown that the evolution of the momentum map from time $0$ to time $T$ is given by the relation
\begin{equation}\label{eq:noether_force}
\left[ \left(\mathbf{J} _L \circ \left(F_L^u\right)^T\right) (q(0),\dot{q}(0))- \mathbf{J} _L(q(0),\dot{q}(0))  \right] \cdot \xi = \int _0^T f_{L}^u(q(t),\dot{q}(t)) \cdot \xi_Q(q(t))\, \text{d} t.
\end{equation}
Equation \eqref{eq:noether_force} shows, that forcing will generally alter the momentum map. However, in the special case that the forcing is orthogonal to the group action, the above relation shows that Noether's theorem will still hold. 

\begin{theorem}[Forced Noether's theorem]
Consider a Lagrangian system $L:TQ \rightarrow \mathbb{R}$ with control forcing $f_{L}: TQ\times U \rightarrow T^*Q$ and a symmetry action $\phi: G \times Q \rightarrow Q$ such that $\langle f_{L}^u(q,\dot{q}), \xi_Q(q)  \rangle = 0$ for all $(q,\dot{q}) \in TQ$, $u\in \mathcal{C}(U)$ and all $\xi \in \mathfrak{g}$. Then the Lagrangian momentum map $\mathbf{J} _L:TQ \rightarrow \mathfrak{g}^*$ will be preserved by the flow, such that $\mathbf{J} _L \circ \left(F_L^u\right)^t= \mathbf{J} _L$ for all $t$.
\end{theorem}

\subsection{Discrete Mechanics}\label{discrete_forcing}

\paragraph{The Discrete Lagrangian.}

Again we consider a configuration manifold $Q$, and define the (``discrete'') state space to be $Q \times Q$. Rather than considering a configuration $q$ and velocity $\dot{q}$ (or momentum $p$), we now consider two configurations $q_0$ and $q_1$, which should be thought of as two points on a curve $q$ which are a time step $h>0$ apart, i.e.\ $q_0 \approx q(0)$ and $q_1\approx q(h)$. The manifold $Q \times Q$ is locally isomorphic to $TQ$ and thus contains the same amount of information.

A \emph{discrete Lagrangian} is a function $L_d: Q \times Q \rightarrow \mathbb{R}$, which we think of as approximating the action integral along the exact solution curve segment $q$ between $q_0$ and $q_1$:
\begin{equation}\label{discreteLagrangian}
L _{d}(q_0,q_1) \approx \int _0^h L ( q (t), \dot{q} (t) ) \,dt.
\end{equation}
We consider the grid $\{t_k = kh\,|\, k=0,\dots,N\}$, $Nh=T$, and define the {\em discrete path space} 
\[
\mathcal{C}_d(Q) = \mathcal{C}_d(\{ t_k\}_{k=0}^N,Q) = \{ q_d: \{ t_k\}_{k=0}^N \rightarrow Q\}.
\]
We will identify a discrete trajectory $q_d\in \mathcal{C}_d(Q)$ with its image $q_d = \{ q_k\}_{k=0}^N$, where $q_k = q_d(t_k)$. The \emph{discrete action map} $\mathfrak{G}_d: \mathcal{C}_d(Q)\rightarrow \mathbb{R}$ along this sequence is calculated by summing the discrete Lagrangian on each adjacent pair and defined by 
\[\mathfrak{G}_d (q_d) = \sum\limits_{k=0}^{N-1} L_d(q_k,q_{k+1}).\]
As the discrete path space $\mathcal{C}_d$ is isomorphic to $Q\times \cdots \times Q$ ($N+1$ copies), it can be given a smooth product manifold structure. The discrete action $\mathfrak{G}_d$ inherits the smoothness of the discrete Lagrangian $L_d$.

The tangent space $T_{q_d}\mathcal{C}_d(Q)$ to $\mathcal{C}_d(Q)$ at $q_d$ is the set of maps $v_{q_d}: \{ t_k\}_{k=0}^N \rightarrow TQ$ such that $\tau_q \circ v_{q_d} = q_d$, which we will denote by $v_{q_d} = \{ (q_k,v_k) \}_{k=0}^N$.

To complete the discrete setting for forced mechanical systems, we present a discrete formulation of the control forces introduced in the previous section. Since the control path $u:[0,T]\rightarrow U$ has no geometric interpretation, we have to find an appropriate discrete formulation to identify a discrete structure for the Lagrangian control force.

\paragraph{Discrete Lagrangian Control Forces.}\label{discLangcf}

Analogous to the replacement of the path space by a discrete path space, we replace the control path space by a discrete one. To this end we consider a refined  grid $\Delta \tilde{t}$, generated via a set of control points $0\leq c_1<\cdots < c_s\leq 1$ as $\Delta \tilde{t}=\{t_{k\ell}=t_k+c_\ell h \,|\, k=0,\ldots,N-1,\ell=1,\ldots,s\}$. With this notation the \index{discrete control path space|textit}{\em discrete control path space} is defined to be
\[
\mathcal{C}_d(U) = \mathcal{C}_d(\Delta\tilde{t},U) = \{ u_d: \Delta \tilde{t} \rightarrow U\}.
\]
 We define the {\em intermediate control samples} $u_k$ on $[t_k,t_{k+1}]$ as $u_k=(u_{k1},\dots,u_{ks}) \in U^s$ to be the values of the control parameters guiding the system from $q_k = q_d(t_k)$ to $q_{k+1}= q_d(t_{k+1})$, where $u_{kl} = u_d(t_{kl})$ for $l\in \{1,\dots,s\}$.




With this definition of the discrete control path space, we take two \index{discrete Lagrangian control force|textit}{\em discrete Lagrangian control forces} $f_{d}^+, f_{d}^-: Q \times Q \times U^s \rightarrow T^*Q$, given in coordinates as
\begin{align}
&f_{d}^+(q_k,q_{k+1},u_k) = (q_{k+1},f_{d}^+(q_k,q_{k+1},u_k)),\\
&f_{d}^-(q_k,q_{k+1},u_k) = (q_k,f_{d}^-(q_k,q_{k+1},u_k)),
\end{align}
also called \index{left and right discrete forces|textit}{\em left and right discrete forces}.\footnote{Observe that the discrete control force is now dependent on the discrete control path.} Analogously to the continuous case, we interpret the two discrete Lagrangian control forces as two families of discrete fiber-preserving Lagrangian forces $ f_d^{u_k,\pm}: Q\times Q \rightarrow T^*Q$ in the sense that $\pi_Q\circ f_d^{u_k,\pm } = \pi^\pm_Q$ with fixed $u_k\in U^s$ and with the projection operators $\pi^+_Q: Q\times Q\rightarrow Q, (q_k,q_{k+1})\mapsto q_{k+1}$ and $\pi^-_Q: Q\times Q\rightarrow Q, (q_k,q_{k+1})\mapsto q_k$.
We combine the two discrete control forces to give a single one-form $f_d^{u_k}: Q\times Q \rightarrow T^*(Q\times Q)$ defined by
\begin{equation}
f_{d}^{u_k}(q_k,q_{k+1})\cdot (\delta q_k, \delta q_{k+1}) = f_{d}^{u_k,+}(q_k,q_{k+1})\cdot \delta q_{k+1} + f_{d}^{u_k,-}(q_k,q_{k+1}) \cdot \delta q_k,
\end{equation}
where $f_{d}(q_k,q_{k+1},u_k)$ denotes the family of all one-forms $f_d^{u_k}(q_k,q_{k+1})$ with fixed $u_k\in U^s$.
To simplify the notation we denote the left and right discrete forces by $f_k^\pm :=  f_{d}^{\pm}(q_k,q_{k+1},u_k)$, respectively, and the pair consisting of both by $f_k := f_{d}(q_k,q_{k+1},u_k)$.

\begin{wrapfigure}[8]{r}[0in]{0.6\linewidth}
\vspace{-5ex}
\begin{center}
\includegraphics[scale=0.45,angle=0]{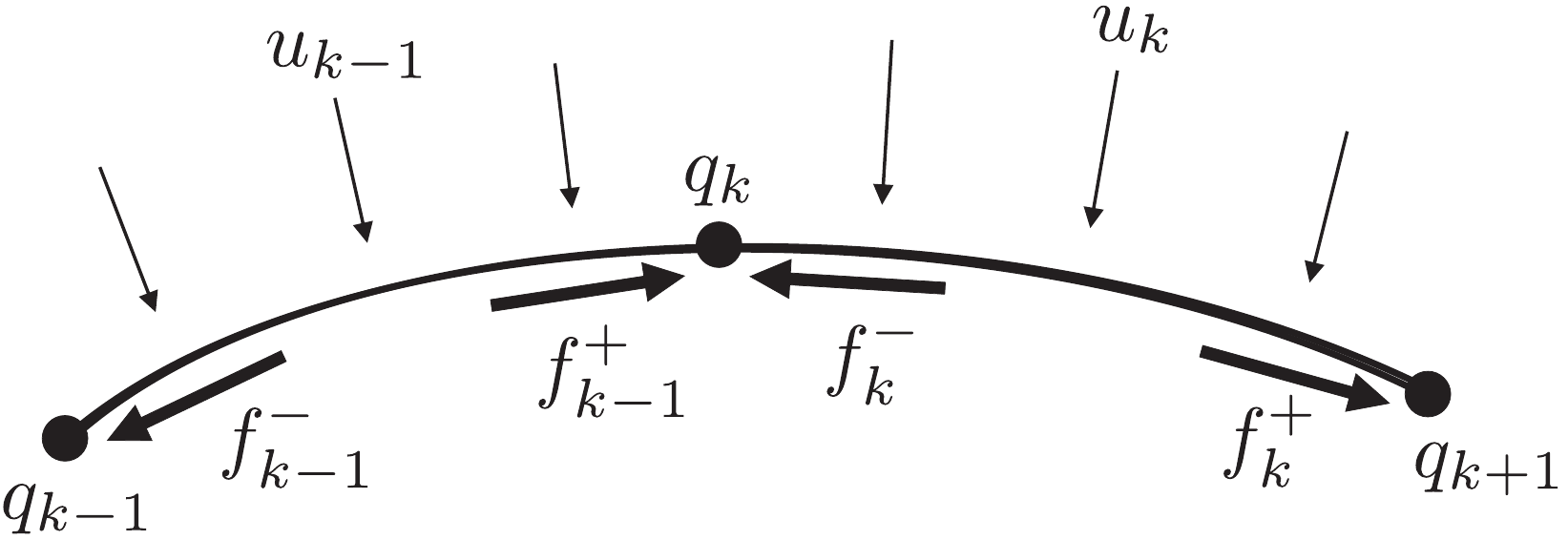}
\caption[Left and right discrete forces]{Left and right discrete forces.}
\label{fig:discreteforces}
\end{center}
\end{wrapfigure}

 Referring to Figure \ref{fig:discreteforces}, we interpret the left discrete force $f_{k-1}^+$ as the force resulting from the continuous control force acting during the time span $[t_{k-1},t_{k}]$ on the configuration node $q_k$. The right discrete force $f_k^-$ is the force acting on $q_k$ resulting from the continuous control force during the time span $[t_k,t_{k+1}]$.
 
\paragraph{The Discrete Lagrange-d'Alembert Principle.}
As with discrete Lagrangians, the discrete control forces also depend on the time step $h$, which is important when relating discrete and continuous mechanics. Given such forces, we modify the discrete Hamilton's principle, following \cite{Kane00}, to the {\em discrete Lagrange-d'Alembert principle}, which seeks discrete curves $\{q_k\}_{k=0}^N$ that satisfy
\begin{equation}\label{eq:varmech_dldap}
\delta \sum\limits_{k=0}^{N-1} L_d(q_k,q_{k+1}) + \sum\limits_{k=0}^{N-1} \left[ f_{d}^-(q_k,q_{k+1},u_k)\cdot \delta q_k +  f_{d}^+(q_k,q_{k+1},u_k)\cdot \delta q_{k+1} \right] = 0
\end{equation}
for all variations $\{\delta q_k\}_{k=0}^N$ vanishing at the endpoints. This is equivalent to the {\em forced discrete Euler-Lagrange equations}
\begin{equation}\label{eq:varmech_FDEL}
D_2 L_d (q_{k-1},q_k) + D_1 L_d (q_{k},q_{k+1}) + f_{d}^+(q_{k-1},q_k,u_{k-1}) + f_{d}^-(q_{k},q_{k+1},u_{k}) = 0.
\end{equation}
These equations implicitly define the forced discrete Lagrangian map $F_{L_d}^{u_{k-1},u_k}: Q\times Q  \rightarrow Q\times Q$ for fixed controls $u_{k-1}, u_k \in U^s$,
mapping $(q_{k-1},q_k)$ to $(q_k,q_{k+1})$.

The discrete Lagrangian one-forms $\Theta^+_{L_d}$ and  $\Theta^-_{L_d}$ are in coordinates
\begin{subequations}
\begin{align}
\Theta^+_{L_d}(q_0,q_1)& = D_2 L_d(q_{0},q_1) \mathbf{d}  q_1 = \frac{\partial L_d}{\partial q_1^i} dq_1^i,\\
\Theta^-_{L_d}(q_0,q_1)& = -D_1 L_d(q_{0},q_1) \mathbf{d}  q_0 = -\frac{\partial L_d}{\partial q_0^i} dq_0^i. 
\end{align}
\end{subequations}

In the absence of external forces, the discrete Lagrangian maps inherit the properties of symplectic preservation from the continuous Lagrangian flows. That means the discrete Lagrangian symplectic form $\Omega_{L_d} = \mathbf{d}  \Theta_{L_d}^+ =  \mathbf{d}  \Theta_{L_d}^-$,  with coordinate expression
\[\Omega_{L_d}(q_0,q_1) = \frac{\partial^2L_d}{\partial q_0^i \partial q_1^j} \mathbf{d}  q_0^i \wedge \mathbf{d}  q_1^j \] 
is preserved under the discrete Lagrangian map as $(F_{L_d})^* (\Omega_{L_d}) = \Omega_{L_d}$, if no external forcing is present.

\paragraph{The Discrete Legendre Transforms with Forces.}
Although in the continuous case we used the standard Legendre transform for systems with forcing, in the discrete case it is necessary to take the {\em forced discrete Legendre transforms} 
\begin{subequations}\label{eq:varmech_dlegf}
\begin{align}
&\mathbb{F}^{f+} L_d: (q_0,q_1,u_0) \mapsto (q_1,p_1) = (q_1, D_2 L_d(q_0,q_1) + f_{d}^+(q_0,q_1,u_0)),\label{eq:varmech_dlegf+}\\
&\mathbb{F}^{f-} L_d: (q_0,q_1,u_0) \mapsto (q_0,p_0) = (q_0, -D_1 L_d(q_0,q_1)-f_{d}^-(q_0,q_1,u_0)).\label{eq:varmech_dlegf-}
\end{align}
\end{subequations}
Again, we denote with $\mathbb{F}^{f\pm} L_d^{u_0}$ the forced discreteLegendre transforms for fixed controls $u_0 \in U^s$.
Using these definitions and the forced discrete Euler-Lagrange equations \eqref{eq:varmech_FDEL}, we can see that the corresponding \index{forced discrete Hamiltonian map|textit}{\em forced discrete Hamiltonian map} $$\tilde{F}_{L_d}^{u_0} = \mathbb{F}^{f\pm} L_d^{u_1} \circ F_{L_d}^{u_0,u_1} \circ (\mathbb{F}^{f\pm} L_d^{u_0})^{-1}$$ is given by the map $\tilde{F}_{L_d}^{u_0}: (q_0,p_0) \mapsto (q_1,p_1)$, where
\begin{subequations}\label{eq:varmech_impf}
\begin{align}
& p_0 = -D_1L_d(q_0,q_1) - f_{d}^{u_0,-}(q_0,q_1),\label{eq:varmech_impf-}\\
& p_1 = D_2 L_d(q_0,q_1) + f_{d}^{u_0,+}(q_0,q_1)\label{eq:varmech_impf+},
\end{align}
\end{subequations}
which is the same as the standard discrete Hamiltonian map with the discrete forces added.

\begin{figure}[htbp]
\begin{center}
\begin{tikzpicture}[fill=blue!20]

\path (2.5,4) node(a) {$(q_0,q_1)$};
\path (7.5,4) node(b)  {$(q_1,q_2)$};
\path (0,0) node(c) {$(q_0,p_0)$};
\path (5,0) node(d) {$(q_1,p_1)$};
\path (10,0) node(e) {$(q_2,p_2)$};

\path (5.0,4.0) node[anchor=south] (f) {$F_{L_d}^{u_0,u_1}$};
\path (1.25,2) node[anchor=east] (g) {$\mathbb{F}^{f-}L_d^{u_0}$};
\path (3.9,2) node[anchor=west] (h) {$\mathbb{F}^{f+}L_d^{u_0}$};
\path (6.3,2) node[anchor=west] (i) {$\mathbb{F}^{f-}L_d^{u_1}$};
\path (8.9,2) node[anchor=west] (j) {$\mathbb{F}^{f+}L_d^{u_1}$};
\path (2.5,0) node[anchor=north] (k) {$\tilde{F}_{L_d}^{u_0}$};
\path (7.5,0) node[anchor=north] (l) {$\tilde{F}_{L_d}^{u_1}$};
\draw[thick,black,|->] (a) -- (b);
\draw[thick,black,|->] (a)  -- (c);
\draw[thick,black,|->] (a)  -- (d);
\draw[thick,black,|->] (b)  -- (d);
\draw[thick,black,|->] (b)  -- (e);
\draw[thick,black,|->] (c)  -- (d);
\draw[thick,black,|->] (d)  -- (e);

\end{tikzpicture}
\caption[Correspondence between the forced discrete Lagrangian and the forced discrete Hamiltonian map]{Correspondence between the forced discrete Lagrangian and the forced discrete Hamiltonian map.}\label{diagram:LH}
\label{fig:discretemaps}
\end{center}
\end{figure}
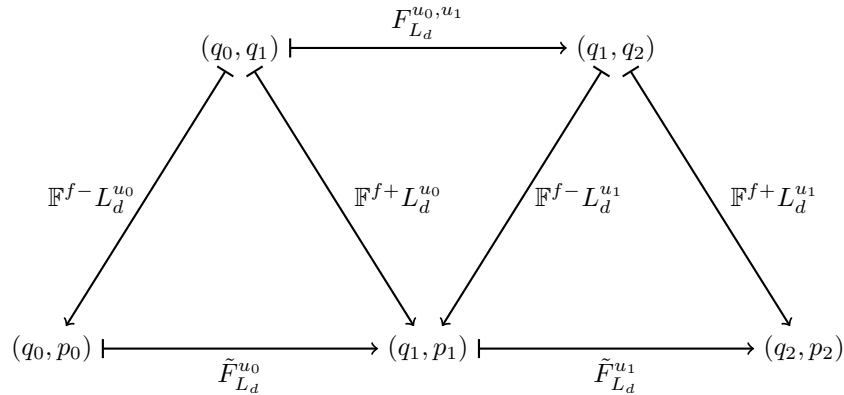

Figure~\ref{fig:discretemaps} shows that the following two definitions of the forced discrete Hamiltonian map
\begin{subequations}
\begin{align}
\tilde{F}_{L_d}^{u_0} & = \mathbb{F}^{f\pm} L_d^{u_1} \circ F_{L_d}^{u_0,u_1} \circ (\mathbb{F}^{f\pm}L_d^{u_0})^{-1},\\
\tilde{F}_{L_d}^{u_0} & = \mathbb{F}^{f+}L_d^{u_0}\circ (\mathbb{F}^{f-}L_d^{u_0})^{-1},\label{eq:def_discHam2}
\end{align}
\end{subequations}
are equivalent with coordinate expression \eqref{eq:varmech_impf}. Thus from expression \eqref{eq:def_discHam2} and Figure~\ref{diagram:LH}, it becomes clear, that the forced discrete Hamiltonian map that maps $(q_0,p_0)$ to $(q_1,p_1)$ depends on $u_0$ only.

\paragraph{The Discrete Noether Theorem with Forcing.}
As in the unforced case, we can formulate a discrete version of the forced Noether's theorem (for the derivation see for example \cite{MaWe01}). To this end, the discrete momentum map in presence of forcing is defined as
\begin{align*}
\mathbf{J} _{L_d}^{f+}(q_0,q_1) \cdot \xi &=  \langle \mathbb{F}^{f+} L_d^{u_0}(q_0,q_1), \xi_Q(q_1) \rangle,\\
\mathbf{J} _{L_d}^{f-}(q_0,q_1) \cdot \xi &=  \langle \mathbb{F}^{f-} L_d^{u_0}(q_0,q_1), \xi_Q(q_0) \rangle.
\end{align*} 
The evolution of the discrete momentum map is described by
\begin{equation}\label{eq:noether_force_d}
\left[ \mathbf{J} _{L_d}^{f+} \circ \left(F_{L_d}^{u_d} \right)^{N-1} - \mathbf{J} _{L_d}^{f-}\right] (q_0,q_1) \cdot \xi = \sum\limits_{k=0}^{N-1} f_{d}^{u_k}(q_k,q_{k+1}) \cdot \xi_{Q\times Q}(q_k,q_{k+1}).
\end{equation}
Again, in the case that the forcing is orthogonal to the group action we have the unique momentum map $\mathbf{J} _{L_d}^{f}: Q\times Q \rightarrow \mathfrak{g}^*$ and it holds:

\begin{theorem}[Forced discrete Noether's theorem]\label{th:discforcNoe}
Consider a discrete Lagrangian system $L_d:Q\times Q \rightarrow \mathbb{R}$ with discrete control forces $f_{d}^+, f_{d}^-: Q\times Q\times U^s \rightarrow T^*Q$ and a symmetry action $\phi: G \times Q \rightarrow Q$ such that $\langle f_{d}^{u_k}, \xi_{Q\times Q}  \rangle = 0$ for all $\xi \in \mathfrak{g}$ and $u_k\in U^s,\, k\in\{0,\ldots,N-1\}$. Then the discrete Lagrangian momentum map $\mathbf{J} ^f_{L_d}:Q \times Q \rightarrow \mathfrak{g}^*$ will be preserved by the discrete Lagrangian evolution map, such that $\mathbf{J} _{L_d}^f \circ F_{L_d}^{u_k,u_{k+1}} = \mathbf{J} _{L_d}^f$.
\end{theorem}

\subsection{The Discrete vs. the Continuous Lagrangian Systems}\label{corres_cd}

In this section, we relate the continuous and the discrete Lagrangian system.
First, along the lines of \cite{MaWe01}, we define expressions for the discrete mechanical objects that exactly reflect the continuous ones. Based on the exact discrete expressions, we determine the order of consistency concerning the difference between the continuous and the discrete mechanical system. 

\paragraph{Exact Discrete Lagrangian and Forcing.}

Given a regular Lagrangian $L:TQ \rightarrow \mathbb{R}$ and a Lagrangian control force $f_{L}:TQ\times U \rightarrow T^*Q$, we define the \index{exact discrete Lagrangian|textit}{\em exact discrete Lagrangian} $L_d^E: Q\times Q \times \mathbb{R} \rightarrow \mathbb{R}$ and the exact discrete control forces $f_{d}^{E+}, f_{d}^{E-}: Q\times Q \times \mathcal{C}([0,h],U) \times \mathbb{R}\rightarrow T^*Q$ to be
\begin{align}
L_d^E(q_0,q_1,h) &= \int _0^h L(q(t),\dot{q}(t)) \, \text{d} t, \label{eq:varmech_exactL}\\
f_{d}^{E+}(q_0,q_1,\mathrm{u}_0,h) &= \int _0^h f_{L}(q(t),\dot{q}(t),u(t)) \cdot \frac{\partial q(t)}{\partial q_1}\, \text{d} t, \label{eq:varmech_exactf+}\\
f_{d}^{E-}(q_0,q_1,\mathrm{u}_0,h) &= \int _0^h f_{L}(q(t),\dot{q}(t),u(t)) \cdot \frac{\partial q(t)}{\partial q_0}\, \text{d} t, \label{eq:varmech_exactf-}
\end{align}
with $\mathrm{u}_k\in \mathcal{C}([kh,(k+1)h],U)$ and where $q:[0,h]\rightarrow Q$ is the solution of the forced Euler-Lagrange equations \eqref{eq:varmech_FEL} with control function $u:[0,h]\rightarrow U$ for $L$ and $f_{L}$ satisfying the boundary conditions $q(0)=q_0$ and $q(h) = q_1$. Observe, that the exact discrete control forces depend on an entire control path in contrast to the continuous control forces. Consequently, the {\em exact forced discrete Legendre transforms} are given by
\begin{align*}
&\mathbb{F}^{f+}L_d^E(q_0,q_1,\mathrm{u}_0,h) = (q_1,D_2 L_d^E(q_0,q_1,h) +  f_{d}^{E+}(q_0,q_1,\mathrm{u}_0,h) ),\\
&\mathbb{F}^{f-}L_d^E(q_0,q_1,\mathrm{u}_0,h) = (q_0,-D_1 L_d^E(q_0,q_1,h) -  f_{d}^{E-}(q_0,q_1,\mathrm{u}_0,h) ).
\end{align*}
As in \S\ref{discrete_forcing} $\mathbb{F}^{f\pm}L_d^{E,\mathrm{u}_k}(q_k,q_{k+1})$ and $f_{d}^{E,\mathrm{u}_k,\pm}(q_k,q_{k+1})$ denote the exact discrete forces and the exact forced discrete Legendre transforms for a fixed $\mathrm{u}_k\in \mathcal{C}([kh,(k+1)h],U)$.
The next lemma is based on a result in \cite{MaWe01} (Lemma 1.6.2) extended to the presence of control forces and establishes a special relationship between the Legendre transforms of a regular Lagrangian and its corresponding exact discrete Lagrangian. This also proves that exact discrete Lagrangians are automatically regular.

\begin{lemma}
A regular Lagrangian $L$ and the corresponding exact discrete Lagrangian $L_d^E$ have Legendre transforms related by
\begin{align*}
&\mathbb{F}^{f+}L_d^E(q_0,q_1,\mathrm{u}_0,h) = \mathbb{F} L(q_{0,1}(h),\dot{q}_{0,1}(h)),\\
&\mathbb{F}^{f-}L_d^E(q_0,q_1,\mathrm{u}_0,h) = \mathbb{F} L(q_{0,1}(0),\dot{q}_{0,1}(0)),
\end{align*}
for sufficiently small $h$ and close $q_0, q_1 \in Q$. Here $q_{0,1}$ denotes the solution of the corresponding Euler-Lagrange equations with $q(0)=q_0, q(h)=q_1$.
\end{lemma}
\begin{proof} Analogous to the proof in \cite{MaWe01} (Lemma 1.6.2) for the unforced case we begin with $\mathbb{F}^{f-}L_d^E$ and compute
\begin{align*}
  \mathbb{F}^{f-}L_d^E(q_0,q_1,\mathrm{u}_0,h)   
&   = -\int _0^h \left[  \frac{\partial L}{\partial q} \cdot \frac{\partial q_{0,1}}{\partial q_0}  + \frac{\partial L}{\partial \dot{q}}  \cdot \frac{\partial \dot{q}_{0,1}}{\partial q_0} \right]\, \text{d} t - \int _0^h f_{L}\cdot \frac{\partial q_{0,1}}{\partial q_0}\, \text{d} t \\[5pt]
&   =  - \int _{0}^h \left[ \frac{\partial L}{\partial q} - \frac{\text{d}}{\text{d} t}  \frac{\partial L}{\partial \dot{q}} - f_{L}  \right] \cdot \frac{\partial q_{0,1}}{\partial q_0}\, \text{d} t - \left[ \frac{\partial L}{\partial \dot{q}}\cdot \frac{\partial q_{0,1}}{\partial q_0}  \right]_0^h,
\end{align*}
using integration by parts. Since $q_{0,1}(t)$ is a solution of the forced Euler-Lagrange equations the first term is zero. With $q_{0,1}(0)= q_0$ and $q_{0,1}(h)= q_1$ for the second term we get 
\[  \frac{\partial q_{0,1}}{\partial q_0}(0)= \mathrm{id} \quad \text{and}\quad  \frac{\partial q_{0,1}}{\partial q_0}(h)=0. \]
Substituting these into the above expression for $\mathbb{F}^{f-}L_d^E$ now gives
\[\mathbb{F}^{f-}L_d^E(q_0,q_1,\mathrm{u}_0,h) = \frac{\partial L}{\partial \dot{q}}(q_{0,1}(0),\dot{q}_{0,1}(0)),\]
which is simply the definition of $\mathbb{F} L(q_{0,1}(0),\dot{q}_{0,1}(0))$.

The result for $\mathbb{F}^{f+}L_d^E$ can be established by a similar computation. \end{proof}

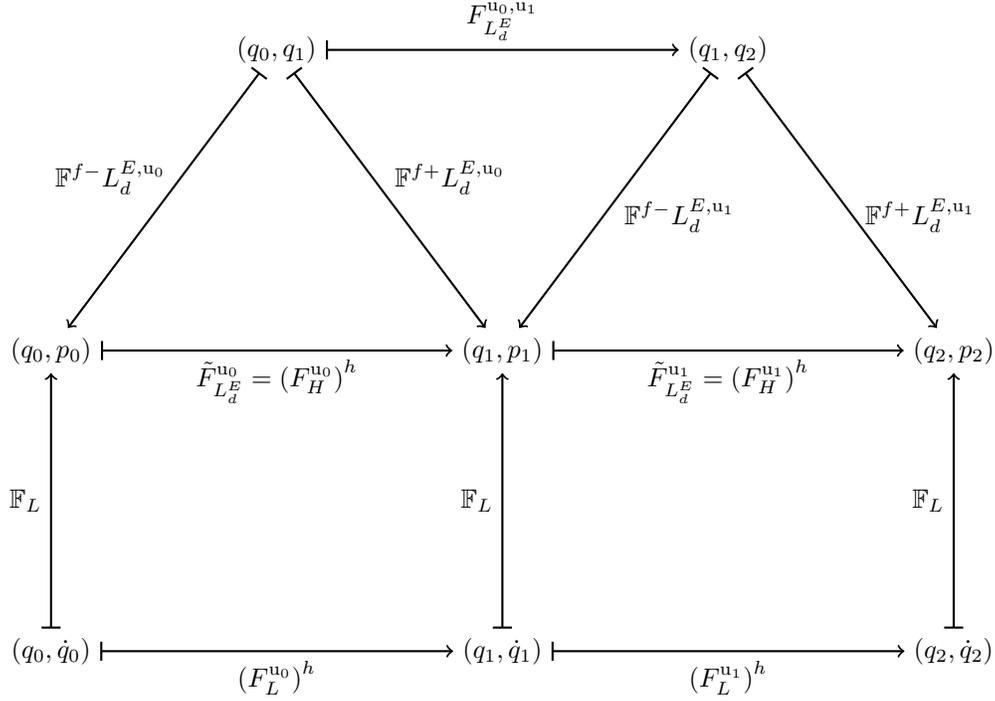
\begin{figure}[h!]
\begin{center}
\begin{tikzpicture}
\path (2.0,4) node(a) {$(q_0,q_1)$};
\path (8,4) node(b)  {$(q_1,q_2)$};
\path (-1,0) node(c) {$(q_0,p_0)$};
\path (5,0) node(d) {$(q_1,p_1)$};
\path (11,0) node(e) {$(q_2,p_2)$};
\path (-1,-4) node(cc) {$(q_0,\dot{q}_0)$};
\path (5,-4) node(dd) {$(q_1,\dot{q}_1)$};
\path (11,-4) node(ee) {$(q_2,\dot{q}_2)$};

\path (5.0,4.0) node[anchor=south] (f) {$F_{L_d^E}^{\mathrm{u}_0,\mathrm{u}_1}$};
\path (0.65,2.3) node[anchor=east] (g) {$\mathbb{F}^{f-}L_d^{E,\mathrm{u}_0}$};
\path (3.45,2.3) node[anchor=west] (h) {$\mathbb{F}^{f+}L_d^{E,\mathrm{u}_0}$};
\path (6.5,1.8) node[anchor=west] (i) {$\mathbb{F}^{f-}L_d^{E,\mathrm{u}_1}$};
\path (9.7,1.8) node[anchor=west] (j) {$\mathbb{F}^{f+}L_d^{E,\mathrm{u}_1}$};
\path (2.0,0) node[anchor=north] (k) {$\tilde{F}_{L_d^E}^{\mathrm{u}_0} = \left(F_H^{\mathrm{u}_0}\right)^h$};
\path (8,0) node[anchor=north] (l) {$\tilde{F}_{L_d^E}^{\mathrm{u}_1} =  \left(F_H^{\mathrm{u}_1}\right)^h$};

\path (-1,-2) node[anchor=east] (k) {$\mathbb{F}_L$};
\path (5,-2) node[anchor=east] (l) {$\mathbb{F}_L$};
\path (11,-2) node[anchor=east] (l) {$\mathbb{F}_L$};

\path (2.0,-4) node[anchor=north] (k) {$\left(F_L^{\mathrm{u}_0}\right)^h$};
\path (8,-4) node[anchor=north] (l) {$\left(F_L^{\mathrm{u}_1}\right)^h$};
\draw[thick,black,|->] (a) -- (b);
\draw[thick,black,|->] (a)  -- (c);
\draw[thick,black,|->] (a)  -- (d);
\draw[thick,black,|->] (b)  -- (d);
\draw[thick,black,|->] (b)  -- (e);
\draw[thick,black,|->] (c)  -- (d);
\draw[thick,black,|->] (d)  -- (e);

\draw[thick,black,|->] (cc)  -- (c);
\draw[thick,black,|->] (cc)  -- (dd);
\draw[thick,black,|->] (dd)  -- (d);
\draw[thick,black,|->] (dd)  -- (ee);
\draw[thick,black,|->] (ee)  -- (e);

\end{tikzpicture}
\caption[Correspondence between the exact discrete Lagrangian and forces and the continuous forced Hamiltonian flow]{Correspondence between the exact discrete Lagrangian and forces and the continuous forced Hamiltonian flow.}\label{diagram:LEH}
\end{center}
\end{figure}

Combining this result with the diagram in Figure \ref{diagram:LH} gives the commutative diagram shown in Figure \ref{diagram:LEH} for the exact discrete Lagrangian and forces.
The diagram also clarifies the following observation, that was already proved in \cite{MaWe01} (Theorem 1.6.3) for unforced systems and can now be established for the forced case as well: Consider the pushforward of both, the continuous Lagrangian and forces and their exact discrete Lagrangian and discrete forces to $T^*Q$, yielding a forced Hamiltonian system with Hamiltonian $H$ and a forced discrete Hamiltonian map $\tilde{F}_{L_d^E}^{\mathrm{u}_k}$, respectively. Then, for a sufficiently small time step $h\in \mathbb{R}$, the forced  Hamiltonian flow map equals the pushforward discrete Lagrangian map: $\left(F_H^{\mathrm{u}_0}\right)^h = \tilde{F}_{L_d^E}^{\mathrm{u}_0}$.

\paragraph{Order of Consistency.} In the previous paragraph we observed that the exact discrete Lagrangian and forces generate a forced discrete Hamiltonian map that exactly equals the forced Hamiltonian flow of the continuous system. 
Since we are interested in using discrete mechanics to reformulate optimal control problems, we generally do not assume that $L_d$ and $L$ or $H$ are related by \eqref{eq:varmech_exactL}. Moreover, the exact discrete Lagrangian and exact discrete forces are generally not computable. In this section we determine the error we obtain by using discrete approximations for the Lagrangian and the control forces.

Rather than considering how closely the trajectory of $F$ matches the exact trajectory given by $F_H$, as is done  in forward error analysis, we use the concept of variational error analysis as introduced in \cite{MaWe01}. In this context we consider how closely a discrete Lagrangian matches the exact discrete Lagrangian given by the action. For forced systems, we additionally take into account how closely the discrete forces match the exact discrete forces. As was stated in the last section, if the discrete Lagrangian is equal to the action and the discrete forces are given by \eqref{eq:varmech_exactf+} and \eqref{eq:varmech_exactf-}, then the corresponding forced discrete Hamiltonian map $\tilde{F}_{L_d}^{u_k}$ will exactly equal the forced flow $F_H^{\mathrm{u}_k}$. Usually, this is just an approximation, therefore we define the local variational error as follows (see \cite{MaWe01}). For the following we assume $Q$ to be a normed vector space equipped with the norm $\|\cdot \|$.

\begin{definition}[Order of consistency]
A given discrete Lagrangian $L_d$ is {\em of order} $r$, if for a fixed curve $u\in\mathcal{C}(U)$ there exist an open subset $V_v\subset TQ$ with compact closure and constants $C_v > 0$ and $h_v>0$ such that
\begin{equation}\label{eq:DMOC_orderL}
| L_d(q(0),q(h),h) - L_d^E(q(0),q(h),h) | \le C_v h^{r+1} 
\end{equation}
for all solutions $q(t)$ of the forced Euler-Lagrange equations with initial condition $(q^0,\dot{q}^0) \in V_v$ and for all $h\le h_v$.
Analogously, a given discrete force $f_{d}^{u_0,\pm}$ is {\em of order} $r$, if there exist an open subset $V_w\subset TQ$
 with compact closure and constants $C_w > 0$ and $h_w>0$ such that
\begin{equation}\label{eq:DMOC_orderf}
\| f_{d}^{u_0,\pm}(q(0),q(h),h) - f_{d}^{E,\mathrm{u}_0,\pm}(q(0),q(h),h) \| \le C_w h^{r+1} 
\end{equation}
for all solutions $q(t)$ of the forced Euler-Lagrange equations with initial condition $(q^0,\dot{q}^0) \in V_w$ for all $h\le h_w$ and with fixed $u_0\in U^s$ and $\mathrm{u}_0\in \mathcal{C}([0,h],U)$. The discrete Legendre transforms $\mathbb{F}^+L_d^{u_0}$ and $\mathbb{F}^-L_d^{u_0}$ of a discrete Lagrangian $L_d$ are {\em of order} $r$ if there exist an open subset $V_f\subset TQ$ with compact closure and constants $C_f>0$ and $h_f>0$ such that
\begin{equation}\label{eq:DMOC_orderLeg}
\| \mathbb{F}^{\pm} L_d^{u_0}(q(0),q(h),h) - \mathbb{F}^{\pm}L_d^{E,\mathrm{u}_0}(q(0),q(h),h) \| \le C_f h^{r+1} 
\end{equation}
for all solutions $q(t)$ of the forced Euler-Lagrange equations with initial condition $(q^0,\dot{q}^0) \in V_f$ and for all $h\le h_f$ and with fixed $u_0\in U^s$ and $\mathrm{u}_0\in \mathcal{C}([0,h],U)$.
\end{definition}

To give a relationship between the orders of a discrete Lagrangian, discrete forces, the forced discrete Legendre transforms, and their forced discrete Hamiltonian maps, we first have to introduce what we mean by {\em equivalence} of discrete Lagrangians: $L_d^1$ is {\em equivalent} to $L_d^2$ if their discrete Hamiltonian maps are equal. For the forced case, we say analogously, that for fixed $u_k\in U^s$ the discrete pair $(L_d^1,f_{d}^{u_k,1})$ is equivalent to the discrete pair $(L_d^2,f_{d}^{u_k,2})$ if their forced discrete Hamiltonian maps are equal, such that $\tilde{F}_{L_d^{1}}^{u_k} = \tilde{F}_{L_d^2}^{u_k}$. With $\tilde{F}_{L_d^1}^{u_k} = \mathbb{F}^{f+}L_d^{u_k,1} \circ (\mathbb{F}^{f-}L_d^{u_k,1})^{-1}$, it follows that if $(L_d^1,f_{d}^{u_k,1})$ and $(L_d^2,f_{d}^{u_k,2})$ are equivalent, then their forced discrete Legendre transforms are equal. Thus, equivalent pairs of discrete Lagrangians and control forces generate the same integrators. 

The following theorem is an extended version of the unforced case stated in \cite{MaWe01} (Theorem 2.3.1).

\begin{theorem} \label{th:ordercalculation}
Given a regular Lagrangian $L$, a Lagrangian control force $f_{L}$ and a corresponding Hamiltonian $H$ with Hamiltonian control force $f_{H}$. Then for a fixed curve $u\in \mathcal{C}(U)$ and a corresponding fixed control sequence $u_k\in U^s$ the following statements are equivalent for a discrete Lagrangian $L_d$ and the discrete forces $f_{d}^{\pm}$:
\begin{itemize}
\item[(i)] the forced discrete Hamiltonian map for $(L_d,f_{d}^{u_k,\pm})$ is of order $r$,
\item[(ii)] the forced discrete Legendre transforms of $(L_d,f_{d}^{u_k,\pm})$ are of order $r$,
\item[(iii)] $(L_d,f_{d}^{u_k,\pm})$ is equivalent to a pair of discrete Lagrangian and discrete forces, both of order $r$.
\end{itemize}
\end{theorem}

\begin{proof} Considering the proof of the unforced version in \cite{MaWe01} (Theorem 2.3.1), this proof is straightforward by taking into account the discrete Lagrangian control forces. For an alternative approach of the proof of Theorem 2.3.1 in \cite{MaWe01} see \cite{PatCu}. \end{proof}



 Note that, given a discrete Lagrangian and discrete forces, their order can be calculated by expanding the expressions for $L_d(q(0),q(h),h)$ and $f_{d}^{u_k,\pm}$ in a Taylor series in $h$ and comparing these to the same expansions for the exact Lagrangian and the exact forces, respectively. If the series agree up to $r$ terms, then the discrete Lagrangian is of order $r$. Analogously, the discrete forces are of order $r$, if for both expansions the first $r$ terms are identical.

\section{Optimal Control of a Mechanical System}

\subsection{The continuous setting}\label{sec:locp}

On the configuration space $Q$ we consider a mechanical system described by a regular Lagrangian $L: TQ \rightarrow \mathbb{R}$. Additionally, assume that a Lagrangian control force acts  on the system and is defined by a map $f_{L}: TQ \times U  \rightarrow T^*Q$ with $f_{L}: (q,\dot{q},u) \mapsto (q,f_{L}(q,\dot{q},u))$ and $u: [0,T]\rightarrow U$, the time-dependent control parameter.  Note that the Lagrangian control force may include both dissipative forces within the mechanical system and external control forces resulting from actuators steering the system.

\paragraph{The Lagrangian Optimal Control Problem.} We now consider the following optimal control problem: During the time interval $[0, T]$, the mechanical system described by the Lagrangian $L$ is to be moved on a curve $q$ from an initial state $(q(0),\dot{q}(0))=(q^{0},\dot q^{0}) \in TQ$ to a final state. The motion is influenced via a Lagrangian control force $f_{L}$ with control parameter $u$ such that a given \emph{objective functional} 
\begin{equation}
\label{eq:DMOC_J}
J(q,u) = \int_0 ^T C(q(t), \dot{q}(t), u(t))\, \text{d} t + \Phi(q(T),\dot{q}(T))
\end{equation}
is minimized. Here $C: TQ \times U \rightarrow \mathbb{R}$ and  $\Phi: TQ \rightarrow \mathbb{R}$ (Mayer term) are continuously differentiable cost functions. The final state $(q(T), \dot{q}(T))$ is required  to fulfil a constraint $r(q(T), \dot{q}(T),q^T,\dot{q}^T)=0$ with $r:TQ \times TQ \rightarrow \mathbb{R}^{n_r}$ and $(q^T,\dot{q}^T) \in TQ$ given.

The motion of the system is to satisfy the Lagrange-d'Alembert principle, which requires that
\begin{equation}\label{eq:DMOC_ldap}
\delta \int_0^T L(q(t),\dot{q}(t)) \, \text{d} t + \int_0^T f_{L}(q(t),\dot{q}(t),u(t))\cdot \delta q(t) \,\text{d} t = 0
\end{equation}
for all variations $\delta q$ with $\delta q(0)=\delta q(T)=0$. In many cases, one encounters additional constraints on the states and controls given by $h(q(t),\dot{q}(t),u(t))\ge 0$ with $h:TQ \times U \rightarrow \mathbb{R}^{n_h}$.

To summarize, we are faced with the following 
\begin{problem}[Lagrangian optimal control problem (LOCP)]\label{prob:locp}
\begin{subequations}\label{eq:LOCP}
\begin{equation}\label{eq:LOCP1}
\min_{q\in \mathcal{C}(Q),u\in\mathcal{C}(U)}  J(q,u)
\end{equation}
subject to
\begin{eqnarray}
\delta \int_0^T L(q(t),\dot{q}(t)) \, \text{d} t + \int_0^T f_{L}(q(t),\dot{q}(t),u(t))\cdot \delta q(t) \, \, d t &=& 0, \\\label{eq:LOCP2}
(q(0),\dot{q}(0)) &=& (q^0,\dot{q}^0), \\ \label{eq:LOCP3}
h(q(t),\dot{q}(t), u(t)) &\ge& 0, \quad t\in [0,T],\label{eq:LOCP4}\\
r(q(T),\dot{q}(T),q^T,\dot{q}^T) &=& 0.\label{eq:LOCP5}
\end{eqnarray} 
\end{subequations} 
\end{problem}
The interval length $T$ may either be fixed, or appear as degree of freedom in the optimization problem.

\begin{definition}
A curve $(q,u)\in \cC(Q)\times\cC(U)$ is {\em admissible} (or \index{feasible|textit}{\em feasible}), if it fulfills the constraints \eqref{eq:LOCP2}--\eqref{eq:LOCP5}. The set of all admissible (feasible) curves is the {\em admissible (feasible) set} of Problem \ref{prob:locp}. An admissible (feasible) curve $(q^*,u^*)$ is an {\em optimal solution} of Problem \ref{prob:locp}, if 
\begin{equation}\label{eq:opt_sol}
J(q^*,u^*) \le J(q,u)
\end{equation}
for all admissible (feasible) curves $(q,u)$. An admissible (feasible) curve $(q^*,u^*)$ is a {\em local optimal solution}, if (\ref{eq:opt_sol}) is fulfilled in a neighborhood of $(q^*,u^*)$.  The function $q^*$ is called {\em (locally) optimal trajectory}, and $u^*$ is the {\em (locally) optimal control}.
\end{definition}

\paragraph{The Hamiltonian Optimal Control Problem.}

We now formulate the problem using the Hamiltonian variant for the system dynamics. This is equivalent to the Lagrangian formulation as we have seen in paragraph \S\ref{varmech} on the Legendre transform with forces. 

For a set $\mathcal{R} = \{(q,\dot{q}) \in TQ \;| \; g(q,\dot{q}) \ge 0  \}$ determined via a constraint $g: TQ\rightarrow \mathbb{R}^{n_g}$ on  $TQ$ we obtain the corresponding set in the cotangent bundle as
$\tilde{\mathcal{R}}= \{(q,p) \in T^*Q \;| \; \tilde g (q,p) \ge 0  \}$ with $\tilde g = g \circ (\mathbb{F} L)^{-1}$. Correspondingly, the optimal control problem in the Hamiltonian formulation reads as follows:

\begin{problem}[Hamiltonian optimal control problem (HOCP)]\label{prob:hocp}
\begin{subequations}\label{eq:HOCP}
\begin{equation}\label{eq:HOCP1}
\min_{(q,p)\in T^* \mathcal{C}(Q),u\in\mathcal{C}(U)}  \tilde{J}(q,p,u) = \int_0 ^T \tilde{C}(q(t), p(t), u(t))\, \, d t  + \tilde{\Phi}(q(T),p(T))
\end{equation}
subject to
\begin{eqnarray}
\dot{q}(t) &=&  \nabla_p H(q(t),p(t)), \label{eq:HOCP2}\\
\dot{p}(t) &=&  -\nabla_q H(q(t),p(t)) + f_{H}(q(t),p(t),u(t)),\\\label{eq:HOCP3}
(q(0),p(0)) &=& (q^0,p^0),\\ \label{eq:HOCP4}
\tilde{h}(q(t),p(t), u(t)) &\geq& 0, \quad t\in [0,T] \\ \label{eq:HOCP5}
\tilde{r}(q(T),p(T),q^T,p^T) &=& 0.\label{eq:HOCP6}
\end{eqnarray}
\end{subequations} 
where $(q^T,p^T)  = \mathbb{F} L(q^T,\dot{q}^T)$, $p(0)  =D_2L(q(0),\dot{q}(0))$, $p^0 =D_2L(q^0,\dot{q}^0)$ and $\tilde{\Phi} = \Phi \circ (\mathbb{F} L)^{-1}$ etc.
\end{problem}

\paragraph{Necessary Optimality Conditions.} 

In this paragraph, we introduce necessary conditions for the optimality of a solution $(q^*,u^*)\in \cC(Q)\times\cC(U)$ to Problems \ref{prob:locp} and \ref{prob:hocp}. We restrict ourselves to the case of problems with the controls pointwise constrained to the (nonempty) set  $U = \{u \in \mathbb{R}^{n_u} \,|\, h(u) \ge 0\} $ and fixed final time $T$. With $x=(q,p)$ and
\[
f(x,u) = (\nabla_p H(q,p),  -\nabla_q H(q,p) + f_{H}(q,p,u))
\]
we can rewrite \eqref{eq:HOCP} as
\begin{subequations}\label{eq:ocp_pont}
\begin{equation}\label{eq:ocp_pont1}
\min_{x\in T^* \cC(Q),u\in\cC(U)} J(x,u) = \int_{0}^T \tilde{C}(x(t),u(t)) \, \text{d} t + \Phi(x(T))
\end{equation}
subject to
\begin{eqnarray}
\dot{x}(t) &=& f(x(t),u(t)), \\
x(0) &=& x_0,\\
u(t) &\in& U, \quad t\in [0,T]\\ \label{eq:ocp_pont2}
\tilde{r}(x(T),x^T) &=& 0.
\end{eqnarray}
\end{subequations}
Necessary conditions for optimality of solution trajectories $\eta(\cdot) = (x(\cdot), u(\cdot))$ can be derived based on variations of an augmented cost function, the {\em Lagrangian} of the system:
\begin{equation}\label{eq:ocp_Jaug}
\mathcal{L} (\eta,\lambda) = \int_{0}^T \tilde{C}(x(t),u(t)) + \lambda^T(t) \cdot (\dot{x} - f(x(t),u(t))) \, \, d t + \Phi(x(T)),
\end{equation}
where the variables $\lambda_i,\; i=1,\dots,n_x$, are the \index{adjoint variables}{\em adjoint variables} or the \index{costates}{\em costates}.
A point $(\eta^*,\lambda^*)$ is a {\em saddle point} of \eqref{eq:ocp_Jaug}, if for all $\eta$ and $\lambda$ it holds
\[\mathcal{L}(\eta,\lambda^*) \le \mathcal{L}(\eta^*,\lambda^*) \le \mathcal{L}(\eta^*,\lambda).
\]
The function 
\begin{equation}
\mathcal{H}(x,u,\lambda) := -\tilde{C}(x,u) + \lambda^T \cdot f(x,u)
\end{equation}
is called the {\em Hamiltonian} of the optimal control problem.
When setting variations of $\mathcal{L}$ with respect to $\eta$ and $\lambda$ to zero, the resulting Euler-Lagrange equations provide \index{necessary optimality condition for optimal control problem}necessary optimality condition for the optimal control problem \eqref{eq:ocp_pont}. Formally, one obtains the following celebrated Theorem (cf.\ \cite{PBGM62}):

\begin{theorem}[Pontryagin Maximum Principle]\label{th:pontryagin}
Let $(x^*,u^*)$ be an optimal solution to \eqref{eq:ocp_pont}. Then there exists a piecewise continuously differentiable function $\lambda: [0,T] \rightarrow \mathbb{R}^{n_x}$ and a vector $\alpha \in \mathbb{R}^{n_r}$ such that 
\begin{subequations}
\begin{equation}
\mathcal{H}(x^*(t),u^*(t),\lambda(t)) = \max\limits_{u\in U} \mathcal{H}(x(t),u,\lambda(t)) \quad t\in[0,T],
\end{equation}
and $\lambda$ solves the following initial value problem:
\begin{eqnarray}
\lambda(T) &=&\nabla_x \Phi(x^*(T)) -\nabla_x \tilde{r}(x^*(T),x^T) \, \alpha,\\
\dot{\lambda} &=& -\nabla_x \mathcal{H}(x^*,u^*,\lambda).
\end{eqnarray}
\end{subequations}
\end{theorem}


\paragraph{Transformation to Mayer Form.} Within the previous sections we introduced optimal control problems in \emph{Bolza form}, in which the objective functional consists of a cost functional of integral form and a final point constraint. For the error analysis in \S\ref{corres_ocp} it will be useful to transform the problem into {\em Mayer form}, in which the objective functional consists of the final point constraint only.  To this end we introduce a new state variable as
\begin{equation}\label{eq:Bolza_toMayer}
z(t) := \int_{0}^t C(q(\tau),\dot{q}(\tau),u(\tau))\, \text{d}\tau,\quad 0\le t\le T,
\end{equation}
resp.\ ${y}(t) := \int_{0}^t \tilde{C}(q(\tau),p(\tau),u(\tau))\, \text{d}\tau$, $0\le t\le T$ for the Hamiltonian form. By extension of the state space from $TQ$ to $TQ\times \mathbb{R}$, and from $T^*Q$ to $T^*Q\times \mathbb{R}$, respectively, the new objective function in Mayer form reads 
\[ 
J(q,\dot{q},u) = z(T) + \Phi(q(T),\dot{q}(T))
\]
resp.\ $\tilde{J}(q,p,u) = {y}(T) + \tilde{\Phi}(q(T),p(T))$. Equation~(\ref{eq:Bolza_toMayer}) is typically adjoint to the problem description as an additional differential equation of the form
\begin{equation}\label{eq:mayer_lang_cont}
\dot{z} = C(q,\dot{q},u), \quad z(0) = 0,
\end{equation}
resp.\ $\dot{{y}} = \tilde{C}(q,p,u)$, $y(0) = 0$.

\subsection{The discrete setting}\label{sec:oc_dm}

For the numerical solution we need a discretized version of Problem \ref{prob:locp}. To this end we formulate an optimal control problem for the discrete mechanical system described by discrete variational mechanics introduced in \S\ref{discrete_forcing}.
In \S\ref{corres_ocp} we show how the optimal control problem for the continuous and the discrete mechanical system are related.

To obtain a discrete formulation, we replace each expression in \eqref{eq:LOCP} by its discrete counterpart in terms of discrete variational mechanics.
As described in \S\ref{discrete_forcing}, we replace the state space $TQ$ of the system by $Q\times Q$ and a path $q:[0,T] \to Q$ by a discrete path $q_{d}:\{0,h,2h,\ldots,Nh=T\}\to Q$ with $q_{k}=q_{d}(kh)$.  Analogously,  the continuous control path $u:[0,T]\to U$ is replaced by a discrete control path $u_{d}: \Delta \tilde{t} \to U$ (writing $u_{k}=(u_{d}(kh+c_\ell h))_{\ell=1}^s \in U^s$).

\paragraph{The Discrete Lagrange-d'Alembert Principle.} 

Based on this discretization, the action integral in \eqref{eq:DMOC_ldap} is approximated on a time slice $[kh,(k+1)h]$ by the \emph{discrete Lagrangian} $L_{d}:Q\times Q\to\mathbb{R}$, as in equation \eqref{discreteLagrangian},
\[
L_{d}(q_{k},q_{k+1}) \approx \int_{kh}^{(k+1)h} L(q(t),\dot q(t))\, \text{d} t,
\]
and likewise the virtual work by the left and right discrete forces,
\[
f_{k}^-\cdot \delta q_{k} + f_{k}^+ \cdot \delta q_{k+1}
\approx \int_{kh}^{(k+1)h} f_{L}(q(t),\dot{q}(t),u(t))\cdot \delta q(t)\, \text{d} t,
\]
where $f_{k}^-, f_{k}^+\in T^*Q$. 

As introduced in equation \eqref{eq:varmech_dldap}, the discrete version of the Lagrange-d'Alembert principle \eqref{eq:DMOC_ldap} requires one to find discrete paths $\{q_k\}_{k=0}^{N}$ such that for all variations $\{\delta q_{k}\}_{k=0}^N$ with $\delta q_{0}=\delta q_{N}=0$, one has
\begin{equation}\label{ddap}
\delta \sum_{k=0}^{N-1} L_d(q_k,q_{k+1})+ \sum_{k=0}^{N-1} \left[f_k^-\cdot \delta q_k + f_k^+ \cdot \delta q_{k+1} \right]= 0,
\end{equation}
or, equivalently, the forced discrete Euler-Lagrange equations
\begin{equation}\label{DEL}
D_2L_d(q_{k-1},q_{k}) + D_1L_d(q_{k},q_{k+1}) +f_{k-1}^+ +  f_k^- = 0,
\end{equation} 
$k=1,\ldots,N-1$. 

\paragraph{Boundary Conditions.}\label{DMOC:bc}
In the next step, we need to incorporate the boundary conditions ${q}(0)={q}^0, \dot{{q}}(0)=\dot{{q}}^0$ and $r({q}(T),\dot{q}(T),q^T,\dot{q}^T)=0$ into the discrete description. Those on the configuration
level can be used as constraints in a straightforward way as
${q}_0={q}^0$. However, 
since in the present
formulation velocities are approximated in a time interval $[t_k,t_{k+1}]$ (as opposed to an
approximation at the time nodes), the velocity
conditions  have to be transformed to conditions on the conjugate momenta. These are defined at each time node using the discrete Legendre transform.
The presence of forces at the time nodes has to be incorporated into that transformation leading to the forced discrete Legendre transforms $\mathbb{F}^{{f}^-}L_d$ and $\mathbb{F}^{f+}L_d$ defined in \eqref{eq:varmech_dlegf}.
Using the standard Legendre transform $\mathbb{F} L:TQ \rightarrow T^*Q, (q,\dot{q}) \mapsto (q,p) = (q,D_2L(q,\dot{q}))$ leads to the discrete initial constraint on the conjugate momentum 
\[D_2L(q^0,\dot{q}^0) + D_1 L_d(q_0,q_1) + f_{d}^-(q_0,q_1,u_0) =0.\]
As shown in the previous section, we can transform the boundary condition from a formulation with configuration and velocity to a formulation with configuration and conjugate momentum. 
Thus, instead of considering a discrete version of the final time constraint $r$ on $TQ$ we use a discrete version of the final time constraint $\tilde{r}$ on $T^*Q$. We define the \index{discrete boundary condition|textit}{\em discrete boundary condition} on the configuration level to be 
\[{r}_d: Q \times Q  \times U^s \times TQ \rightarrow \mathbb{R}^{n_r}, \]
\[  {r}_d(q_{N-1},q_N,u_{N-1},q^T,\dot{q}^T) =  \tilde{r} \left(\mathbb{F}^{f+}L_d(q_{N-1},q_N,u_{N-1}), \mathbb{F} L(q^T,\dot{q}^T)\right),\]
with $(q_N,p_N) = \mathbb{F}^{f+}L_d(q_{N-1},q_N,u_{N-1})$ and $(q^T,p^T) = \mathbb{F} L(q^T,\dot{q}^T)$, that is $p_N = D_2 L_d(q_{N-1},q_N) + f_{d}^+(q_{N-1},q_{N},u_{N-1})$ and $p^T=D_2L(q^T,\dot{q}^T)$.

Notice that for the simple final velocity constraint  $$r(q(T),\dot{q}(T),q^T,\dot{q}^T) = \dot{q}(T) - \dot{q}^T,$$ we obtain for the transformed condition on the momentum level $\tilde{r}(q(T),p(T),q^T,p^T) =  p(T)-p^T$ the discrete constraint
\begin{equation}\label{rem:inifinalvel}
-D_2L(q^T,\dot{q}^T)+D_2L_d(q_{N-1},q_{N}) +f_{d}^+(q_{N-1},q_N,u_{N-1}) = 0.
\end{equation}

\paragraph{Discrete Path Constraints.} Opposed to the final time constraint we approximate the path constraint in \eqref{eq:LOCP4} on each time interval $[t_k,t_{k+1}]$ rather than at each time node. Thus, we maintain the formulation on the velocity level and replace the continuous path constraint $h(q(t),\dot{q}(t),u(t)) \ge 0$ by a {\em discrete path constraint} $h_d: Q \times Q \times U^s \rightarrow \mathbb{R}^{sn_h}$ with
\[
h_d(q_k,q_{k+1},u_k)\ge 0, \quad k=0,\dots,N-1.
\]

\paragraph{Discrete Objective Function.} Similar to the Lagrangian we approximate the objective functional in \eqref{eq:DMOC_J} on the time slice $[kh,(k+1)h]$ by
\[
C_{d}(q_{k},q_{k+1},u_{k}) \approx \int_{kh}^{(k+1)h} C(q(t),\dot q(t),u(t))\, \text{d} t.
\]
Analogously to the final time constraint, we approximate the final condition via a discrete version $\Phi_d: Q\times Q \times U^s\rightarrow \mathbb{R}$ yielding  the \index{discrete objective function|textit}\emph{discrete objective function}
\[
J_d(q_d,u_d)=\sum_{k=0}^{N-1} C_d(q_k,q_{k+1},u_k) + \Phi_d(q_{N-1},q_N,u_{N-1}).
\]

\paragraph{The Discrete Optimal Control Problem.}
In summary, after performing the above discretization steps, one is faced with the following discrete optimal control problem.

\begin{problem}[Discrete Lagrangian optimal control problem]\label{probd:locp}

\begin{subequations}\label{eq:gDLOCP}
\begin{equation}\label{eq:gDLOCP1}
\min_{q_d,u_d} J_d(q_d,u_d)
\end{equation}
subject to
\begin{align}
q_0 &= q^0,\label{eq:gDLOCP2}\\
D_2L(q^0,\dot{q}^0)+D_1L_d(q_{0},q_{1}) + f_0^-& = 0,\label{eq:gDLOCP3}\\
D_2 L_d (q_{k-1},q_k) + D_1 L_d (q_{k},q_{k+1}) + f_{k-1}^+ + f_k^- &= 0,\quad k=1,\dots, N-1,\label{eq:gDLOCP4}\\
h_d(q_k, q_{k+1}, u_k) &\ge 0,\quad k=0,\dots, N-1,\label{eq:gDLOCP5}\\
r_d(q_{N-1},q_{N},u_{N-1},q^T,\dot{q}^T) &= 0.\label{eq:gDLOCP6}
\end{align}
\end{subequations} 
\end{problem}

Recall that the $f_k^\pm$ are dependent on $u_k \in U^s$. To incorporate a free final time $T$ as in the continuous setting, the step size $h$ appears as a degree of freedom within the optimization problem. However, in the following formulations and considerations we restrict ourselves to the case of fixed final time $T$ and thus fixed step size $h$.

\paragraph{Necessary Optimality Conditions.}
The system \eqref{eq:gDLOCP} results in a constrained nonlinear optimization problem, also called a  \textit{nonlinear programming problem}, that is an objective function has to be minimized subject to algebraic equality and inequality constraints. 
Let $\xi$ be the set of parameters introduced by the discretization of an infinite dimensional optimal control problem. Then, the nonlinear programming problem to be solved is
\begin{equation}\label{NLP}
\begin{array}{c}
\min_\xi \varphi(\xi)\\[5pt]
\text{subject to}\quad a(\xi) = 0,\; b(\xi) \ge 0,
\end{array}
\end{equation}
where $\varphi: \mathbb{R}^n \rightarrow \mathbb{R},\; a:\mathbb{R}^n \rightarrow \mathbb{R}^m,$ and $b:\mathbb{R}^n \rightarrow \mathbb{R}^p$ are continuously differentiable.

We briefly summarize some terminology: A {\em feasible point} is a point $\xi \in \mathbb{R}^n$ that satisfies $a(\xi) = 0$ and $b(\xi)\ge 0$. A \index{local minimum}{\em local minimum} of (\ref{NLP}) is a feasible point $\xi^*$ which has the property that $\varphi(\xi^*) \le \varphi(\xi)$ for all feasible points $\xi$ in a neighborhood of $\xi^*$. A \index{strict local minimum}{\em strict local minimum} satisfies $\varphi(\xi^*) < \varphi(\xi)$ for all neighboring feasible points $\xi \ne \xi^*$.
{\em Active inequality constraints} $b^{\text{act}}(\xi)$ at a feasible point $\xi$ are those components $b_j(\xi)$ of $b(\xi)$ with $b_j(\xi)=0$. We subsume the equality constraints and the active inequalities at a point $\xi$ (known as \index{active set}{\em active set}) in a combined vector function of \index{active constraint}{\em active constraints} as
\[\tilde{a}(\xi):= \left(\begin{array}{c} a(\xi)\\ b^{\text{act}}(\xi)\end{array}\right).\]
Note that the active set may be different at different feasible points.
{\em Regular} points are feasible points $\xi$ that satisfy the condition that the Jacobian of the active constraints, $\nabla \tilde{a}(\xi)^T$, has full rank, that means that all rows of  $\nabla \tilde{a}(\xi)^T$ are linearly independent.

To investigate local optimality in the presence of constraints, we introduce the \index{Lagrangian multiplier}{\em Lagrangian multiplier} vectors $\lambda \in \mathbb{R}^m$ and $\mu \in \mathbb{R}^p$, that are also called \index{adjoint variables}{\em adjoint variables}, and we define the \index{Lagrangian}{\em Lagrangian function} $\tilde{\mathcal{L}}$ by
\begin{equation}\label{eq:ocp_L}
\tilde{\mathcal{L}}(\xi,\lambda,\mu) := \varphi(\xi) - \lambda^T a(\xi) - \mu^T b(\xi).
\end{equation}
We now state a variant of the {\em Karush-Kuhn-Tucker} necessary conditions for local optimality of a point $\xi^*$. These conditions have been derived first by Karush in 1939 (\cite{Kar39}) and independently by Kuhn and Tucker in 1951 (\cite{KuTu51}).
For brevity, we restrict our attention to regular points only.

\begin{theorem}[Karush-Kuhn-Tucker conditions (KKT)]
If a regular point $\xi^* \in \mathbb{R}^n$ is a local optimum of the NLP problem (\ref{NLP}), then there exist unique Lagrange multiplier vectors $\lambda^* \in \mathbb{R}^m$ and $\mu^* \in\mathbb{R}^p$ such that the triple $(\xi^*,\lambda^*,\mu^*)$ satisfies the following necessary conditions:
\begin{align*}
\nabla_\xi \tilde{\mathcal{L}}(\xi^*,\lambda^*,\mu^*) & = 0,\\
a(\xi^*) & = 0,\\
b(\xi^*) & \ge 0,\\
\mu^* &\ge 0,\\
\mu_j^* b_j(\xi^*) &= 0,\quad j=1,\dots, p.
\end{align*}
\end{theorem}
\begin{proof} 
See for example \cite{JarStoer}. 
\end{proof}

A triple $(\xi^*,\lambda^*,\mu^*)$ that satisfies the Karush-Kuhn-Tucker conditions is called a {\em Karush-Kuhn-Tucker point (KKT point)}.

\paragraph{Transformation to Mayer Form.}
As in the continuous setting, the discrete optimal control problem \eqref{eq:gDLOCP} can be transformed into a problem in Mayer form, that is the objective function consists of a final condition only. The transformation is performed on the discrete level to keep the Lagrangian structure of the original problem. We introduce new variables  
\begin{align*}
z_0 =  0, \quad z_\ell& =  \sum_{k=0}^{\ell-1} C_d(q_k,q_{k+1},u_k),\quad \ell=1,\dots, N,
\end{align*}
and rephrase the discrete problem \eqref{eq:gDLOCP} into a problem of Mayer type:
\begin{subequations}\label{eq:gDLOCP_mayer}
\begin{equation}
\min_{q_d,u_d} z_{N} + \Phi_d(q_{N-1},q_N,u_{N-1})
\end{equation}
subject to
\begin{align}
q_0 &= q^0,\\
z_0 & = 0,\label{gDLOCP_mayer_z0}\\
D_2L(q^0,\dot{q}^0)+D_1L_d(q_{0},q_{1}) + f_0^-& = 0,\\
D_2 L_d (q_{k-1},q_k) + D_1 L_d (q_{k},q_{k+1}) + f_{k-1}^+ + f_k^- &= 0,\quad k=1,\dots, N-1,\\
h_d(q_k, q_{k+1}, u_k) &\ge 0,\quad k=0,\dots, N-1,\\
r_d(q_{N-1},q_{N},u_{N-1},q^T,\dot{q}^T) &= 0,\\
z_{k+1} - z_k - C_d(q_k,q_{k+1},u_k)& = 0,\quad k=0,\dots, N-1.\label{gDLOCP_mayer_z}
\end{align}
\end{subequations} 
Thus equations (\ref{gDLOCP_mayer_z0}) and (\ref{gDLOCP_mayer_z}) provide the corresponding discretization for the additional equation of motion \eqref{eq:mayer_lang_cont} resulting from the Mayer transformation of the Lagrangian optimal control problem on the continuous level.

\paragraph{Fixed Boundary Conditions.} Consider the special case of a problem with fixed initial and final configuration and velocities and without path constraints: 
\begin{subequations}\label{eq:LOCP_fb}
\begin{equation}\label{eq:LOCP_fb1}
\min_{q\in \mathcal{C}(Q),u\in \mathcal{C}(U)} \int_0 ^T C(q(t), \dot{q}(t), u(t))\, \text{d} t 
\end{equation}
subject to
\begin{equation}
\delta \int_0^T L(q(t),\dot{q}(t)) \, \text{d} t + \int_0^T f_{L}(q(t),\dot{q}(t),u(t))\cdot \delta q(t) \,\text{d} t = 0,\label{eq:LOCP_fb2}
\end{equation}
\begin{equation}
q(0) = q^0,\; \dot{q}(0) = \dot{q}^0,\; q(T) = q^T,\; \dot{q}(T) =\dot{q}^T\label{eq:LOCP_fb3}.
\end{equation}
\end{subequations} 

A straightforward way to derive initial and final constraints for the conjugate momenta rather than for the velocities from the variational principle directly is stated in the following proposition:

\begin{proposition}
With $(q^0,p^0) =  \mathbb{F} L(q^0,\dot{q}^0)$ and $(q^T,p^T) =\mathbb{F} L(q^N,\dot{q}^N)$ equations \eqref{eq:LOCP_fb2} and \eqref{eq:LOCP_fb3} are equivalent to the following principle with free initial and final variation and with augmented Lagrangian
 \begin{align}
\delta \left( \int_0^T L(q(t),\dot{q}(t)) \, \, d t + p^0 (q(0)-q^0) -p^T (q(T)-q^T)  \right) + \int_0^T f_{L}(q(t),\dot{q}(t),u(t))\cdot \delta q(t) \, \, d t = 0\label{eq:prop_vel}.
\end{align}
 \end{proposition}
 
 \begin{proof} Variations of \eqref{eq:LOCP_fb2} with respect to $q$ and zero initial and final variation $\delta q(0) = \delta q(T) = 0$ together with \eqref{eq:LOCP_fb3} yield
 \begin{subequations}
 \begin{align}
 & \frac{\text{d}}{\text{d} t} \frac{\partial}{\partial \dot{q}} L (q(t), \dot{q}(t)) - \frac{\partial}{\partial q} L (q(t), \dot{q}(t)) = f_{L}(q(t),\dot{q}(t),u(t)),\\
& q(0) = q^0,\; \dot{q}(0) = \dot{q}^0,\; q(T) = q^T,\; \dot{q}(T) =\dot{q}^T. \label{eq:bounds_cont}
 \end{align}
 \end{subequations}
On the other hand variations of \eqref{eq:prop_vel} with respect to $q$ and $\lambda=(p^0,p^T)$ with free initial and final variation lead to
 \begin{subequations}
 \begin{align}
 & \frac{\text{d}}{\text{d} t} \frac{\partial}{\partial \dot{q}} L (q(t), \dot{q}(t)) - \frac{\partial}{\partial q} L (q(t), \dot{q}(t)) = f_{L}(q(t),\dot{q}(t),u(t)),\\
 &\left. \frac{\partial}{\partial \dot{q}} L(q(t),\dot{q}(t)) \right|_{t=0} = p^T,\label{eq:bounds_impuls1}\\
  &\left. \frac{\partial}{\partial \dot{q}} L(q(t),\dot{q}(t)) \right|_{t=T} = p^0,\label{eq:bounds_impuls2}\\
  &q(0) = q^0,\;  q(T) = q^T.
 \end{align}
 \end{subequations}
The Legendre transform applied to the velocity boundary equations in \eqref{eq:bounds_cont} gives the corresponding momenta boundary equations \eqref{eq:bounds_impuls1} and \eqref{eq:bounds_impuls2}.
\end{proof}

On the discrete level we derive the optimal control problem for fixed initial and final configurations and velocities in an equivalent way.
Thus, we consider the discrete principle with discrete augmented Lagrangian
 \begin{equation}
 \delta\left( \sum\limits_{k=0}^{N-1} L_d(q_k,q_{k+1}) + p^0 (q_0-q^0) -p^T (q_N-q^T) \hspace{-3pt} \right)+ \sum\limits_{k=0}^{N-1} \left[ f_k^-\cdot \delta q_k +  f_k^+\cdot \delta q_{k+1} \right] = 0,
 \end{equation}
which, with free initial and final variation $\delta q_0$ and $\delta q_N$, respectively, is equivalent to
 
 \begin{subequations}\label{eq:DLOCP_fb}
 \begin{equation}
\delta \sum\limits_{k=0}^{N-1} L_d(q_k,q_{k+1}) + \sum\limits_{k=0}^{N-1} \left[ f_k^-\cdot \delta q_k +  f_k^+\cdot \delta q_{k+1} \right] = 0,\label{eq:DLOCP_fb2}
\end{equation}
\begin{equation}
q_0 = q^0,\; p^0+D_1L_d(q_{0},q_{1}) + f_0^- = 0, \label{eq:DLOCP_fb3}
\end{equation}
\begin{equation}\label{eq:DLOCP_fb4}
q_N = q^T,\; -p^T + D_2L_d(q_{N-1},q_{N}) + f_{N-1}^+ = 0,
\end{equation} 
 \end{subequations}
where the second equations in \eqref{eq:DLOCP_fb3} and \eqref{eq:DLOCP_fb4} are exactly the discrete initial and final velocity constraints derived in the remark containing equation \eqref{rem:inifinalvel} with $p^0 = D_2 L(q^0,\dot{q}^0)$ and $p^T = D_2 L(q^T,\dot{q}^T)$.
 
We note that this derivation of the discrete initial and final conditions directly gives the same formulation that we found before by first transforming the boundary condition on the momentum level $T^*Q$ and then formulating the corresponding discrete constraints on $Q\times Q\times U^s$.     

\subsection{The Discrete vs the Continuous Problem}\label{DMOC:DMOC_direct}\label{corres_ocp}
 
This section gives an interpretation of the discrete problem as an approximation to the continuous one.  In addition, we identify certain structural properties that the discrete problem inherits from the continuous one.  We determine the consistency order of the discrete scheme and establish a result on the convergence of the discrete solution as the step size goes to zero. 

\paragraph{The Place of DMOC amongst Solution Methods for Optimal Control Problems.} 

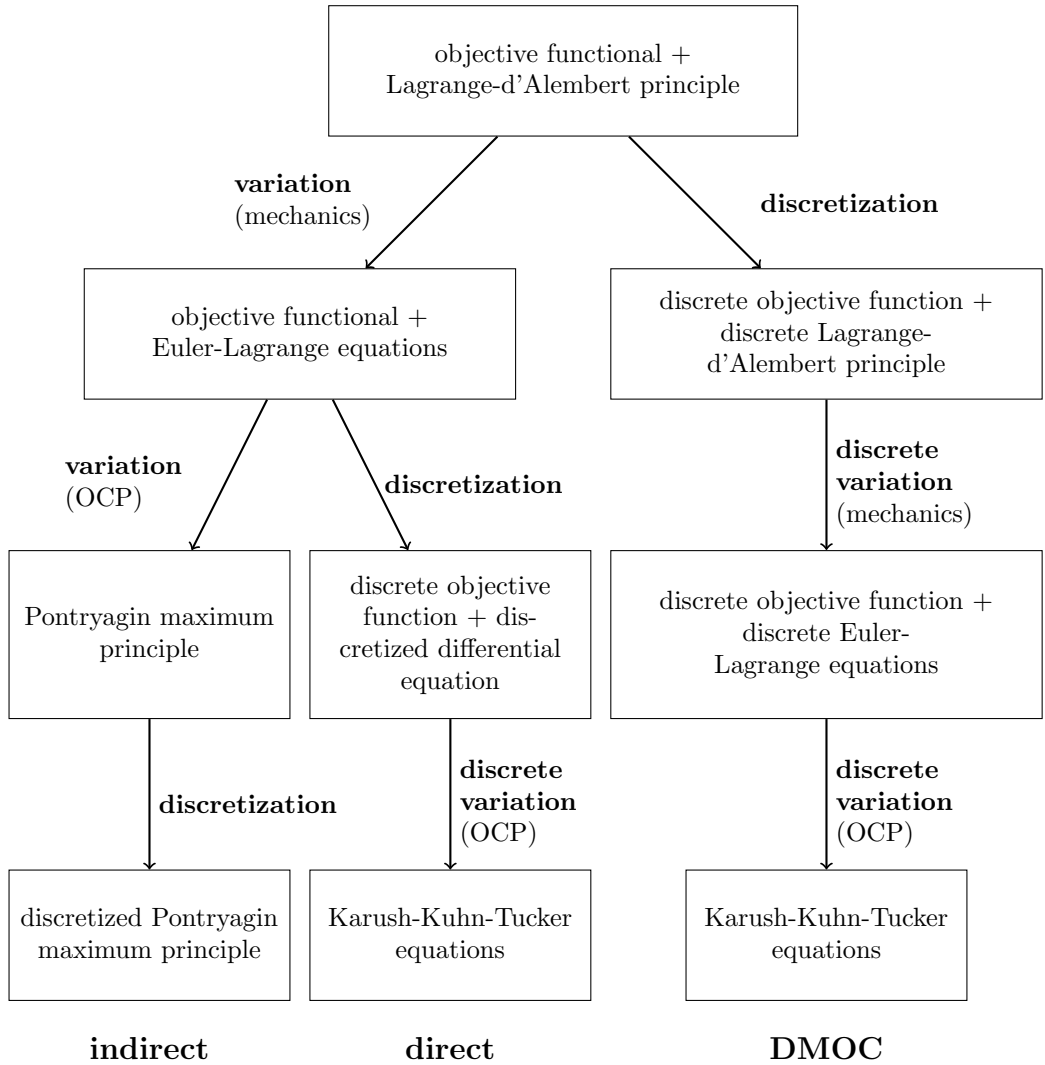
\begin{figure}[h!]
\begin{center}
\begin{tikzpicture}
\path (0.5,4.5) node(a) [rectangle,draw] {\parbox[c][1.5cm][c]{6cm}{\begin{center}objective functional +\\Lagrange-d'Alembert principle\end{center}}};
\path (-1.75,2.75) node[anchor=east] (e) {\parbox{2cm}{{\bf variation}\\(mechanics)}};
\path (-3,1) node(b) [rectangle,draw] {\parbox[c][1.5cm][c]{5.5cm}{\begin{center}objective functional + \\Euler-Lagrange equations\end{center}}};
\path (-2,-1) node[anchor=west] (f) {{\bf discretization}};
\path (-4,-1) node[anchor=east] (f) {\parbox{2cm}{{\bf variation}\\ (OCP)}};
\path (-1,-3) node(c) [rectangle,draw] {\parbox[c][2cm][c]{3.5cm}{\begin{center}discrete objective function + dis-\\cretized differential equation\end{center}}};
\path (-5,-3) node(j) [rectangle,draw] {\parbox[c][2cm][c]{3.5cm}{\begin{center}Pontryagin maximum principle\end{center}}};
\path (-1,-7) node(k) [rectangle,draw] {\parbox[c][1.5cm][c]{3.5cm}{\begin{center}Karush-Kuhn-Tucker equations\end{center}}};
\path (-5,-7) node(l) [rectangle,draw] {\parbox[c][1.5cm][c]{3.5cm}{\begin{center}discretized Pontryagin maximum principle\end{center}}};
\path (3,2.75) node[anchor=west] (g) {{\bf discretization}};
\path (4,1) node(d) [rectangle,draw] {\parbox[c][1.5cm][c]{5.5cm}{\begin{center}discrete objective function +\\ discrete Lagrange-\\d'Alembert principle\end{center}}};
\path (4,-3) node(i) [rectangle,draw] {\parbox[c][2cm][c]{5.5cm}{\begin{center}discrete objective function +\\ discrete Euler-\\Lagrange equations\end{center}}};
\path (4,-7) node(m) [rectangle,draw] {\parbox[c][1.5cm][c]{3.5cm}{\begin{center}Karush-Kuhn-Tucker equations\end{center}}};
\path (4,-1) node[anchor=west] (h) {\parbox{2.75cm}{{\bf discrete\\ variation}\\ (mechanics)}};
\path (4,-5.25) node[anchor=west] (h) {\parbox{2.75cm}{{\bf discrete\\ variation} \\(OCP)}};
\path (-5,-5.25) node[anchor=west] (n) {{\bf discretization}};
\path (-1,-5.25) node[anchor=west] (o) {\parbox{2.75cm}{{\bf discrete\\ variation}\\ (OCP)}};
\path (4,-8.5) node(x)  {\parbox[c][1.5cm][c]{3.5cm}{\begin{center}{\bf {\large DMOC}}\end{center}}};
\path (-5,-8.5) node(x)  {\parbox[c][1.5cm][c]{3.5cm}{\begin{center}{\bf {\large indirect}}\end{center}}};
\path (-1,-8.5) node(x)  {\parbox[c][1.5cm][c]{3.5cm}{\begin{center}{\bf {\large direct}}\end{center}}};
\draw[thick,->] (a) --(b);
\draw[thick,->] (b)  -- (c);
\draw[thick,->] (b)  -- (j);
\draw[thick,->] (d)  -- (i);
\draw[thick,->] (i)  -- (m);
\draw[thick,->] (c)  -- (k);
\draw[thick,->] (j)  -- (l);
\draw[thick,->] (a) --(d);
\end{tikzpicture}
\caption[Optimal control for mechanical systems: the order of variation and discretization for deriving the necessary optimality conditions]{Optimal control for mechanical systems: the order of variation and discretization for deriving the necessary optimality conditions.}\label{fig:intro_comparison}
\end{center}
\end{figure}

In Figure \ref{fig:intro_comparison} we present schematically different discretization strategies for optimal control problems:

\begin{itemize}
\item In an \emph{indirect} method, starting with an objective function and the Lagrange-d'Alembert principle we obtain via two variations (the first for the derivation of the Euler-Lagrange equations and the second for the derivation of the necessary optimality conditions) the Pontryagin maximum principle. The resulting boundary value problem is then solved numerically, e.g. by gradient methods (\cite{Ca1847,Ke60,Tol75,BrHo75,Mi80,ChLu82}), multiple shooting (\cite{Fox60,Kel68,Bu71,Deu74,Bock78,Hilt90}) or collocation (\cite{DiWe75, Bar83, AMR88}). 
 
\item In a \emph{direct} approach, starting form the Euler-Lagrange equations we directly transform the problem into a restricted finite dimensional optimization problem by discretizing the differential equation. Common methods like e.~g.~shooting \cite{Kraft85}, multiple shooting \cite{BoPl84}, or collocation methods \cite{St93}, rely on a direct integration of the associated ordinary differential equations or on 
its fulfillment at certain grid points (see also \cite{Betts98, Pyt99} for an overview of the current state of the art).
The resulting finite dimensional nonlinear constrained 
optimization problem can be solved by standard nonlinear optimization techniques like sequential quadratic 
programming (\cite{Han76, Pow78}).  Implementations are found in software packages like DIRCOL \cite{Str2000}, SOCS \cite{BeHu98}, or MUSCOD \cite{Lei99}).

\item In the \emph{DMOC} approach, rather than discretizing the differential equations arising from the Lagrange-d'Alembert principle, we discretize in the earliest stage, namely already on the level of the variational principle. Then, we consider variations only on the discrete level to derive the restricted optimization problem and its necessary optimality conditions.

This approach derived via the concept of discrete mechanics leads to a special discretization of the system equations based on variational integrators, which are dealt with in detail in \cite{MaWe01}. Thus, the discrete optimal control problem inherits special properties exhibited by variational integrators. In the following, we specify particular important properties and phenomena of variational integrators and try to translate their meaning into the optimal control context.
 
\end{itemize}

\paragraph{Preservation of Momentum Maps.}
If the discrete system, obtained by applying variational integration to a mechanical system, inherits the same symmetry groups as the continuous system, the corresponding discrete momentum maps are preserved. For the forced case the same statement holds, if the forcing is orthogonal to the group action (see Theorem \ref{th:discforcNoe}). 

On the one hand, this means for the optimal control problem, that if the control force is orthogonal to the group action, our discretization leads to a discrete system, for which the corresponding momentum map is preserved. On the other hand, in the case of the forcing not being orthogonal to the group action, the forced discrete Noether's theorem provides an exact coherence between the change in angular momentum and the applied control force via 
\[[\mathbf{J} _{L_d}^{f+} \circ \left(F_{L_d}^{u_d}\right)^{N-1} - \mathbf{J} _{L_d}^{f-}](q_0,q_1) \cdot \xi = \sum\limits_{k=0}^{N-1} f_{d}^{u_k}(q_k,q_{k+1})\cdot \xi_{Q\times Q}(q_k,q_{k+1}),\] 
(see \S\ref{appl} for examples).

\paragraph{Conservation of a Modified Energy.}
Variational integrators are symplectic, which implies that a certain modified energy is conserved (see for example \cite{HaLuWa}). This is an important property if the long time behavior of dynamical systems is considered. For the case of the optimal control of systems with long maneuver time such as low thrust space missions, it would therefore be interesting to investigate the relation between a modified energy and the virtual work. However, this has not been considered within this paper.

\paragraph{Implementation.}
Rather than using a configuration-momentum implementation of variational integrators as proposed in \cite{MaWe01}, we stay on $Q\times Q$. That means we just determine the optimal trajectory for the configuration and the control forces and reconstruct the corresponding momenta and velocities via the  forced discrete Legendre transforms. This yields computational savings. A more detailed description of the computational savings compared to standard discretizations for optimal control problems is given in Remark \ref{rem:savings}.

\subsection{The Correspondence with Runge-Kutta Discretizations}\label{symppartRK}

In this section we are going to show that the discretization derived via the discrete Lagrange-d'Alembert principle is equivalent to one resulting from a finite difference discretization of the associated Hamiltonian system via a symplectic partitioned Runge-Kutta scheme.

\paragraph{Symplectic Partitioned Runge-Kutta Methods.}

As shown in \cite{MaWe01} (Theorem 2.6.1), the discrete Hamiltonian map generated by the discrete Lagrangian is a symplectic partitioned Runge-Kutta method. As we will show, a similar statement is true for discrete Hamiltonian maps with forces. The resulting method is still a partitioned Runge-Kutta method, but no longer symplectic in the original sense since the symplectic form is not preserved anymore due to the presence of control forces. However, we still denote it as a symplectic method having in mind that the symplectic form is preserved only in absence of external control forces.

A partitioned Runge-Kutta method for the regular forced Lagrangian system $(L,f_{L})$ is a map $T^*Q \times U^s \rightarrow T^*Q$ specified by coefficients $b_i, a_{ij}, \tilde{b}_i, \tilde{a}_{ij}$, $i=1,\dots,s$, and defined by $(q_0,p_0,u_0) \mapsto (q_1,p_1)$, where
\parbox{0.2\textwidth}{
\begin{eqnarray*}
q_1 &=& q_0 + h\sum\limits_{j=1}^s b_j \dot{Q}_j, \\  Q_i &= &q_0 + h\sum\limits_{j=1}^s a_{ij} \dot{Q}_j, \\  P_i& =& \frac{\partial L}{\partial \dot{q}}(Q_i,\dot{Q}_i),\end{eqnarray*}} \hfill
\parbox{0.7\textwidth}{
\begin{subequations}\label{eq:DMOC_RK}
\begin{eqnarray}
p_1 &=& p_0 + h\sum\limits_{j=1}^s \tilde{b}_j \dot{P}_j,\label{eq:DMOC_RK1}\\
P_i& =& p_0 + h\sum\limits_{j=1}^s \tilde{a}_{ij} \dot{P}_j, \label{eq:DMOC_RK2}\\
\dot{P}_i &=& \frac{\partial L}{\partial {q}}(Q_i,\dot{Q}_i) + f_{L}(Q_i,\dot{Q}_i,U_i), \label{eq:DMOC_RK3}
\end{eqnarray}
\end{subequations}}
$i=1,\dots,s$, where the points $(Q_i,P_i)$ are known as the \emph{internal stages} and $U_i$ are the control samples given by $U_i = u_{0i}= u_d(t_0+c_ih)$. For $a_{ij} = \tilde{a}_{ij}$ and $b_i = \tilde{b}_i$ the partitioned Runge-Kutta method is a Runge-Kutta method.

The method is symplectic (that is, it preserves the canonical symplectic form $\Omega$ on $T^*Q$ in absence of external forces) if the coefficients satisfy

\hfill\parbox{0.1\textwidth}{
\begin{eqnarray*}\label{con_sym}
 b_i \tilde{a}_{ij} + \tilde{b}_j a_{ji} & = & b_i \tilde{b}_j,\\
 b_i& = & \tilde{b}_i,
 \end{eqnarray*}} 
\parbox{0.7\textwidth}{
\begin{subequations}\label{eq:DMOC_RK_symp}
\begin{eqnarray}
i,j &=& 1,\dots,s,\label{eq:DMOC_RK_symp1}\\
i& =& 1,\dots,s.\label{eq:DMOC_RK2_symp2}
\end{eqnarray}
\end{subequations}}\\
Since discrete Lagrangian maps are symplectic, we can assume that we have coefficients satisfying \eqref{eq:DMOC_RK_symp} and write a discrete Lagrangian and discrete Lagrangian control forces that generate the corresponding symplectic partitioned Runge-Kutta method. Given points $(q_0,q_1) \in Q\times Q$, we can regard \eqref{eq:DMOC_RK} as implicitly defining $p_0, p_1, Q_i, P_i, \dot{Q}_i$ and $\dot{P}_i$ for $i=1,\dots,s$. Taking these to be defined as functions of $(q_0,q_1)$, we construct a discrete Lagrangian
\begin{equation}\label{eq:DMOC_RK_L}
L_d(q_0,q_1,h) = h\sum\limits_{i=1}^s b_i L(Q_i,\dot{Q}_i),
\end{equation}
and left and right discrete forces as
\begin{subequations}\label{eq:DMOC_RK_f}
\begin{align}
f_{d}^+(q_0,q_1,u_0,h) & = h \sum\limits_{i=1}^s b_i f_{L}(Q_i,\dot{Q}_i,U_i) \cdot \frac{\partial Q_i}{\partial q_1}\label{eq:DMOC_RK_f+},\\
f_{d}^-(q_0,q_1,u_0,h) & = h \sum\limits_{i=1}^s b_i f_{L}(Q_i,\dot{Q}_i,U_i) \cdot \frac{\partial Q_i}{\partial q_0}.\label{eq:DMOC_RK_f-}
\end{align}
\end{subequations}
For fixed $u_0$ the corresponding forced discrete Hamiltonian map is exactly the map $(q_0,p_0) \mapsto (q_1,p_1)$ which is the symplectic partitioned Runge-Kutta method for the forced Hamiltonian system \eqref{eq:varmech_forcHam}.

\begin{theorem}\label{equ_L_RK}\label{theroem:RK}
The discrete Hamiltonian map generated by the discrete Lagrangian \eqref{eq:DMOC_RK_L} together with the discrete forces \eqref{eq:DMOC_RK_f} is a partitioned Runge-Kutta method (which is symplectic in the unforced case).
\end{theorem}

\begin{proof} This is straight forward to the proof for unforced systems in  \cite{MaWe01} (Theorem 2.6.1). \end{proof}

\paragraph{Optimal Control and Runge-Kutta Discretizations.}
In this paragraph, we carry forward the results of the previous section into the context of optimal control problems. To formulate the optimal control problem for the discrete system in terms of Runge-Kutta discretizations, we have to give appropriate expressions for the boundary conditions, the path constraints, and the objective functional.
Concerning the dynamics of the mechanical system, we have already seen that the discretization obtained via the discrete Lagrange-d'Alembert principle can be rewritten as a Runge-Kutta scheme for the corresponding mechanical system in terms of $(p,q)$ as in \eqref{eq:DMOC_RK}. Since we consider regular Lagrangians and regular Hamiltonians, we reformulate \eqref{eq:DMOC_RK3} with the help of the Hamiltonian as
\begin{align*}
\dot{Q}_{ki} & = \frac{\partial H}{\partial {p}}(Q_{ki},P_{ki}),\\
\dot{P}_{ki}& = - \frac{\partial H}{\partial {q}}(Q_{ki},P_{ki}) + f_{H}(Q_{ki},P_{ki},U_{ki}),
\end{align*}
where the additional index $k$ denotes the dependence of the intermediate state and control variables on the time interval $[t_k,t_{k+1}].$

\paragraph{Boundary Conditions.} Due to a formulation of the optimal control problem within the Hamiltonian framework, we use the same formulation for the boundary constraint as the one for the continuous Hamiltonian optimal control problem \ref{prob:hocp} evaluated at $(q_N,p_N)$, which reads

\[ 
\tilde{r}(q_N,p_N,q^T,p^T) = 0.
\]

\paragraph{Discrete Path Constraints.} Again we use the formulation of the Hamiltonian optimal control problem \ref{prob:hocp} and enforce the path constraints to hold in each time step $(Q_i,P_i,U_i)$ as

\[
\tilde{h}(Q_{ki},P_{ki},U_{ki}) \ge 0,\quad k=0,\dots,N-1,\quad i=1,\dots,s.
\]

\paragraph{Discrete Objective Function.} We construct the discrete objective function in the same way as the discrete Lagrangian \eqref{eq:DMOC_RK_L}. However, corresponding to the Hamiltonian formulation evaluated in each substep $(Q_i,P_i,U_i)$, it now reads

\[
\tilde{J}_d(q_d,p_d,u_d) = \sum\limits_{k=0}^{N-1} h \sum_{i=1}^s b_i \tilde{C}(Q_{ki},P_{ki},U_{ki}) + \tilde{\Phi}(q_N,p_N),  
\]
where the final constraint holds for the node corresponding to the final time $T$.

Combining terms, the discrete optimal control problem from the Hamiltonian point of view reads
\begin{subequations}\label{eq:dHOCP}
\begin{equation}\label{eq:dHOCP1}
\min_{q_d,p_d,u_d}  \tilde{J}_d(q_d,p_d,u_d) 
\end{equation}
subject to\\
\begin{align}
&q_{k+1} = q_k + h\sum\limits_{j=1}^s b_j \dot{Q}_{kj}, \quad Q_{ki} = q_k + h\sum\limits_{j=1}^s a_{ij} \dot{Q}_{kj}, \quad q_0= q^0,\label{eq:dHOCP2} \\
&p_{k+1} = p_k + h\sum\limits_{j=1}^s \tilde{b}_j \dot{P}_{kj},\quad
P_{ki} = p_k + h\sum\limits_{j=1}^s \tilde{a}_{ij} \dot{P}_{kj}, \quad p_0=p^0, \label{eq:dHOCP3}\\
&\dot{Q}_{ki} = \frac{\partial H}{\partial {p}}(Q_{ki},P_{ki}), \quad \dot{P}_{ki} = - \frac{\partial H}{\partial {q}}(Q_{ki},P_{ki}) + f_{H}(Q_{ki},P_{ki},U_{ki}), \label{eq:dHOCP4}\\
&0 \le \tilde{h}(Q_{ki},P_{ki}, U_{ki}),\quad k=0,\dots,N-1, \;i=1,\dots,s,\label{eq:dHOCP5}\\
&0=\tilde{r}(q_N,p_N,q^T,p^T).\label{eq:dHOCP6}
\end{align}
\end{subequations} 
Problem \eqref{eq:dHOCP} is the finite dimensional optimization problem resulting from the Hamiltonian optimal control problem \ref{prob:hocp} that is discretized via a symplectic partitioned Runge-Kutta scheme.


\begin{theorem}[Equivalence]\label{th:equi}
Given a discrete Lagrangian and discrete forces as defined in (\ref{eq:DMOC_RK_L}) and (\ref{eq:DMOC_RK_f}), respectively. Then, the discrete Lagrangian optimal control problem defined in (\ref{eq:gDLOCP}) and the problem (\ref{eq:dHOCP}) resulting from discretizing the Hamiltonian optimal control problem by a symplectic partitioned Runge-Kutta scheme are equivalent in the sense that both problems have the same set of solutions.
\end{theorem}

\begin{proof}
Assume $(q_d^*,u_d^*)$ is a solution of (\ref{eq:gDLOCP}). By using the discrete Legendre transform we obtain an optimal solution in terms of the discrete momenta as $(q_d^*, p_d^*,u_d^*)$. To prove that $(q_d^*, p_d^*,u_d^*)$ is also an optimal solution for problem (\ref{eq:dHOCP}) we have to show feasibility and optimality.
Theorem \ref{theroem:RK} and application of the forced discrete Legendre transform to the boundary conditions give us the equivalence of equations (\ref{eq:gDLOCP2})--(\ref{eq:gDLOCP4}),(\ref{eq:gDLOCP6}) and (\ref{eq:dHOCP2})--(\ref{eq:dHOCP4}), (\ref{eq:dHOCP6}). For the choice of the discretizations of the curves $q(t)$ and $u(t)$ the inequality condition (\ref{eq:gDLOCP5}) reads as $h(Q_{ki},Q_{k,i+1},U_{ki})\ge 0$ for $i=1,\ldots,s-1$ and $h(Q_{ks},Q_{k+1,0},U_{ks})\ge 0$. Again, due to the forced discrete Legendre transform $h$ and $\tilde{h}$ defined in (\ref{eq:dHOCP5}) determine the same solution set. This makes $(q_d^*, p_d^*,u_d^*)$ a feasible solution of the problem (\ref{eq:dHOCP}). Optimality holds due to $ \tilde{J}_d(q_d^*,p_d^*,u_d^*) =  J_d(q_d^*,u_d^*) \le  J_d(q_d,u_d) = \tilde{J}_d(q_d,p_d,u_d)$ for all feasible points $(q_d,p_d)\in T^* \mathcal{C}_d(Q)$ and $u_d \in \mathcal{C}_d(U)$ (global optimality), or for all feasible points $(q_d,p_d,u_d)$ in a neighborhood of $(q_d^*,p_d^*,u_d^*)$ (local optimality). Analogous, we can show that an optimal solution $(q_d^*,p_d^*,u_d^*)$ of problem (\ref{eq:dHOCP}) is also an optimal solution for problem (\ref{eq:gDLOCP}).
\end{proof}

\paragraph{Mayer Formulation.} Similar to what we did for the discrete Lagrangian optimal control problem, \eqref{eq:dHOCP} can be transformed into an optimal control problem of Mayer type as follows: Analogous to what we did in equation \eqref{eq:gDLOCP_mayer}, we introduce additional variables $y_d = \{y_l\}_{l=0}^{N}$ as
\begin{eqnarray*}
y_0& =& 0,\\
y_l& =& \sum_{k=0}^{l-1} h \sum_{i=1}^s b_i \tilde{C}(Q_{ki},P_{ki},U_{ki}),\quad l=1,\dots, N,
\end{eqnarray*}
yielding the discrete optimal control problem of Mayer type as
\begin{subequations}\label{eq:dHOCP_mayer}
\begin{equation}\label{eq:dHOCP1_mayer}
\min_{q_d,p_d,u_d}  \bar{\Phi}(q_N,p_N,y_{N}) = y_{N}  + \tilde{\Phi}(q_N,p_N)
\end{equation}
subject to (\ref{eq:dHOCP2}), (\ref{eq:dHOCP3}), (\ref{eq:dHOCP4}), (\ref{eq:dHOCP5}) and (\ref{eq:dHOCP6}).
\end{subequations} 
We obtain exactly the same problem by discretizing the continuous Hamiltonian optimal control problem of Mayer type  with the same partitioned Runge-Kutta discretization for the extended system of differential equations
\begin{align*}
\dot{\tilde{q}}(t) &= \nu(\tilde{q}(t),\tilde{p}(t),u(t)),\quad \tilde{q}(0) = \tilde{q}^0,\\
\dot{\tilde{p}}(t) &= \eta(\tilde{q}(t),\tilde{p}(t),u(t)),\quad \tilde{p}(0) = \tilde{p}^0,
\end{align*}
with $\tilde{q} = (q,y)$, $\tilde{p} = p$, $\nu(\tilde{q}(t),\tilde{p}(t),u(t)) = (\nabla_p H(q(t),p(t)), \tilde{C}(q(t),p(t),u(t)))$, $\eta(\tilde{q}(t),\tilde{p}(t),u(t)) = \linebreak -\nabla_p H(q(t),p(t)) + f_{H}(q(t),p(t),u(t))$, $\tilde{q}^0 = (q^0,0)$ and $\tilde{p}^0 = p^0$.

\begin{exams}\label{ex:RK} \quad
\begin{itemize}
\item[{\rm (a)}] {\rm With $b=1$ and $a=\frac{1}{2}$ we obtain the implicit midpoint rule as a symplectic Runge-Kutta scheme, that is the partitioned scheme reduces to a standard one-stage Runge-Kutta scheme. The resulting discretization is equivalent to the discretization derived via the Lagrangian approach with midpoint quadrature.}
\item[{\rm (b)}] {\rm The standard Lobatto IIIA-IIIB partitioned Runge-Kutta method is obtained by using the Lobatto quadrature for the discrete Lagrangian.}
\end{itemize}
\end{exams}

\begin{remark}\label{rem:savings} \quad
\begin{itemize}
\item[{\rm 1.}] {\rm The formulations of the discrete optimal control problem on the one hand based on the discrete Lagrange-d'Alembert principle and on the other hand based on a Runge-Kutta discretization of the Hamiltonian dynamics are equivalent in the sense that the same discrete solution set is described. Note however, that the Lagrangian formulation needs less discrete variables and less constraints for the optimization problem, as it is formulated on the configuration space only, rather than on the space consisting of configurations and momenta: For $q\in\mathbb{R}^n$ and $N$ intervals of discretization we obtain with the Lagrangian approach $(N s + 1)n$ unknown configurations $q_k^\nu,\; k=0,\dots,N-1,\; \nu=0,\dots,s$ with $q_k^s = q_{k+1}^0,\; k=1,\dots,N-1$ and $n  N  (s-1) + n (N-1)$ extended discrete Euler-Lagrange equations, so altogether, $(N s -1) n$ constraints for the optimization problem excluding boundary conditions. The Runge-Kutta approach for the Hamiltonian system yields $2  ( N s +1)n$ unknown configurations and momenta and $2 n  N + 2 n N (s-1)= 2 n  N s$ equality constraints. Thus, via the Runge-Kutta approach we obtain twice as many state variables and $(N s +1)n$ more equality constraints such that the resulting optimization problem is numerically more expensive. Comparisons concerning the computational effort are presented in Section \ref{imp_twolink}.}
\item[{\rm 2.}]
{\rm Another advantage of the variational approach over the sympletic Runge-Kutta discretization is that one does not have to handle any conditions on the coefficients (such as condition (\ref{con_sym})) to enforce symplecticity. The symplecticity of the DMOC discretization results naturally from the variational structure of this approach.} 
\end{itemize}
\end{remark}

\subsection{The Adjoint Systems}\label{adjoints}

The adjoint system provides necessary optimality conditions for a given optimal control problem. In the continuous case the adjoint system is derived via the Pontryagin maximum principle. The KKT equations provide the adjoint system for the discrete optimal control problem. \cite{Hager00} derives a transformed adjoint system using standard Runge-Kutta discretizations. He identifies order conditions on the coefficients of the adjoint scheme up to order $4$. \cite{Bo06} extend these up to order $7$.

By using the same strategy as Hager in \cite{Hager00}, we here show that DMOC leads to a discretization of the same order for the adjoint system as for the state system. Therefore no additional order conditions on the coefficients are necessary.
In the following we ignore path and control constraints and restrict ourselves to the case of unconstrained optimal control problems. 

\paragraph{Continuous Setting.} 

We consider a Hamiltonian optimal control problem in Mayer form  without final constraint, path and control constraints:

\begin{problem}\label{prob_conv_ocp}
\begin{subequations}\label{eq:conv_ocp}
\begin{equation}
\min_{q,p,u}  \Phi(q(T), p(T)) 
\end{equation}\label{eq:conv_ocp1}
subject to
\begin{align}
&\dot{q} =\nu(q,p), \label{eq:conv_ocp2}\\
&\dot{p} = \eta(q,p,u), \label{eq:conv_ocp3}\\ 
& (q(0),p(0)) = (q^0,p^0), \label{eq:conv_ocp4}
\end{align}  
\end{subequations}
with $q,p\in W^{1,\infty}([0,T],\mathbb{R}^n$), $u\in L^\infty([0,T],\mathbb{R}^m)$, $\nu(q,p) = \nabla_p H(q,p)$, $\eta(q,p,u) = -\nabla_q H(q,p) + f_{H}(q,p,u)$.
\end{problem}

 
We also write $x(t) = (q(t),p(t))$, $x^0 = (q^0,p^0)$ and 
$$
\tilde{f}(x,u) = \left(\begin{array}{c} \nu(q,p)\\ \eta(q,p,u) \end{array}\right),
$$ 
such that \eqref{eq:conv_ocp2}-- \eqref{eq:conv_ocp4} reads as
$\dot{x} = \tilde{f}(x,u), x(0) = x^0$.

Along the lines of \cite{Hager00}, we now formulate the assumptions that are employed in the analysis of DMOC discretizations of Problem \ref{prob_conv_ocp}. First, a smoothness assumptions is required to ensure regularity of the solution and the problem functions. Second, we enforce a growth condition that allows for having a unique solution for the control function of the optimal control problem.

\begin{assumption}[Smoothness]\label{ass:smooth}
For some integer $\kappa \ge 2$, Problem \ref{prob_conv_ocp} has a local solution $(x^*,u^*)$ which lies in $W^{\kappa,\infty}\times W^{\kappa-1,\infty}$. There exists an open set $\Omega \subset \mathbb{R}^{2n} \times \mathbb{R}^m$ and $\rho > 0$ such that\footnote{$B_\rho(z)$ is the closed ball centered at $z$ with radius $\rho$.} $B_\rho(x^*(t),u^*(t)) \subset \Omega$ for every $t\in [0,T]$. The first $\kappa$ derivatives of $\nu$ and $\eta$ are Lipschitz continuous in $\Omega$, and the first $\kappa$ derivatives of $\Phi$ are Lipschitz continuous in $B_\rho(x^*(T))$.
\end{assumption}

\paragraph{First order optimality conditions.} 

Under Assumption \ref{ass:smooth} there exits an associated Lagrange multiplier $\psi^* = (\psi^{q,*}, \psi^{p,*})\in W^{\kappa,\infty}$ for which the following form of the first-order optimality conditions derived via the Pontryagin maximum principle is satisfied at $(x^*,\psi^*,u^*)$:
\begin{subequations}\label{eq:DMOC_adjoint}
\begin{align}
&\dot{q} =\nu(q,p),\quad q(0) = q^0,\label{eq:DMOC_adjoint1}\\
&\dot{p} = \eta(q,p,u),\quad p(0) = p^0,\label{eq:DMOC_adjoint2}\\
&\dot{\psi}^q = - \nabla_q \mathcal{H}(q,p,\psi^q,\psi^p,u),\label{eq:DMOC_adjoint3} \quad \psi^q(T) = \nabla_q \Phi(q(T),p(T))\\
&\dot{\psi}^p = - \nabla_p \mathcal{H}(q,p,\psi^q,\psi^p,u),\label{eq:DMOC_adjoint4} \quad \psi^p(T) = \nabla_p \Phi(q(T),p(T))\\
&\nabla_u \mathcal{H}(q(t),p(t),\psi^q(t),\psi^p(t),u(t)) = 0\; \text{for all}\; t\in [0,T],\label{eq:DMOC_adjoint5}
\end{align}
\end{subequations}
with the Hamiltonian $\mathcal{H}$ defined by
\begin{equation}\label{eq:DMOC_Hamiltonian}
\mathcal{H}(q,p,\psi^q,\psi^p,u)  =  \psi^q \nu(q,p) + \psi^p \eta(q,p,u),
\end{equation}
where $\psi^q$ and $\psi^p$ are row vectors in $\mathbb{R}^n$. 
We also write $\mathcal{H}(x,\psi,u) :=  \psi \tilde{f}(x,u)$
such that \eqref{eq:DMOC_adjoint3} and \eqref{eq:DMOC_adjoint4} read as
\[
\dot{\psi} = -\nabla_x \mathcal{H}(x,\psi,u), \quad \psi(T) = \nabla_x\Phi(x(T)).
\]

\paragraph{Second  order optimality conditions.} 

In order to formulate the second-order sufficient optimality conditions we define the following matrices (cf.\ \cite{Hager00}):
\begin{align*}
 A(t) & = \nabla_x \tilde{f}(x^*(t),u^*(t)),\\
 B(t) & = \nabla_u \tilde{f}(x^*(t),u^*(t)),\\
 V(t) & = \nabla^2 \Phi(x^*(T)),\\
 P(t) & = \nabla_{xx} \mathcal{H}(x^*(t),\psi^*(t),u^*(t)),\\
 R(t) & = \nabla_{uu} \mathcal{H}(x^*(t),\psi^*(t),u^*(t)),\\
 S(t) & = \nabla_{xu} \mathcal{H}(x^*(t),\psi^*(t),u^*(t)). 
\end{align*} 
Let $\mathcal{B}$ be the quadratic form defined by
\[\mathcal{B}(x,u) = \frac{1}{2} \left( x(T)^T V x(T) + \langle x, Px \rangle + \langle u, Ru \rangle + 2  \langle x, Su \rangle \right),\]
where  $\langle \cdot,\cdot  \rangle$ denotes the usual $L^2$ inner product. 

\begin{assumption}[Coercivity]
There exists a constant $\alpha > 0$ such that
\[\mathcal{B}(x,u) \ge \alpha \| u\| ^2_{L_2} \quad  \text{for all} \quad (x,u) \in \mathcal{M},\]
where $\mathcal{M} = \left\{ (x,u)\in H^1\times L^2 \mid x(0) = 0 \text{ and } \dot{x} = Ax + Bu \text{ for a.e. } t\in [0,T] \right\}$.
\end{assumption}


Under the assumptions Smoothness and Coercivity, the {\em control uniqueness property} holds (cf.~\cite{Hager79}): the Hamiltonian $\mathcal{H}$ has a locally unique minimizer $u^*=u(x,\psi)$ in the control, depending Lipschitz continuously on $x$ and $\psi$.
Let $\phi = (\phi^q, \phi^p)$ denote the function defined by
\begin{align*}
\phi^q(x,\psi)& = -\nabla_q \mathcal{H}(x,\psi,u^*),\\
\phi^p(x,\psi)& = -\nabla_p \mathcal{H}(x,\psi,u^*).
\end{align*}
Additionally, let $\eta(x,\psi)=\eta(x,u^*)$. By substituting $u^*=u(x,\psi)$ into \eqref{eq:DMOC_adjoint} one obtains a two-point boundary-value problem
\begin{subequations}\label{eq:contbound}
\begin{align}
\dot{q} & = \nu(q,p),\quad q(0)=q^0,\\
\dot{p} & = \eta(q,p,\psi^q,\psi^p),\quad p(0)=p^0,\\
\dot{\psi}^q & = \phi^q(q,p, \psi^q,\psi^p),\quad \psi^q(T) = \nabla_q \Phi(q(T),p(T)),\\
\dot{\psi}^p & = \phi^p(q,p, \psi^q,\psi^p),\quad \psi^p(T) = \nabla_p \Phi(q(T),p(T)).
\end{align}
\end{subequations}

\paragraph{Discrete Setting.}




It was shown in \S\ref{symppartRK} that the discrete Lagrangian control problem is equivalent to Problem \eqref{eq:dHOCP} (resp.\ \eqref{eq:dHOCP_mayer}). 
%
%
Using some transformation and change of variables (see \cite{Hager00} and \cite{OB08}), by reversing the order of time and applying the control uniqueness property, we obtain the following version of the first-order necessary optimality  conditions associated with \eqref{eq:dHOCP} (the Karush-Kuhn-Tucker conditions):

\begin{subequations}\label{eq:stateco}
\begin{align}
q_{k+1} & = q_k + h \sum\limits_{i=1}^s b_i \nu(Q_{ki},P_{ki}), \quad q_0 = q^0, \label{eq:stateco1}\\
p_{k+1} & = p_k + h \sum\limits_{i=1}^s b_i \eta(Q_{ki},P_{ki},\Psi_{ki}^q,\Psi_{ki}^p), \quad {p}_0 = {p}^0, \label{eq:stateco2}\\
Q_{ki} & = q_k + h \sum\limits_{j=1}^s a^q_{ij} \nu(Q_{kj},P_{kj}), \label{eq:stateco3}\\
P_{ki} & = p_k + h \sum\limits_{j=1}^s  a^p_{ij} \eta(Q_{kj},P_{kj},\Psi_{kj}^q, \Psi_{kj}^p), \label{eq:stateco4}\\
\psi_{k+1}^q & =   \psi_{k}^q + h \sum\limits_{i=1}^s b_i  \phi^q(Q_{ki},P_{ki},\Psi_{ki}^q,\Psi_{ki}^p),\quad \psi^q_N = \nabla_q\Phi(q_N,p_N),\label{eq:stateco5}\\
\psi_{k+1}^p & =\psi_{k}^p +  h \sum\limits_{i=1}^s b_i  \phi^p(Q_{ki},P_{ki},\Psi_{ki}^q,\Psi_{ki}^p),\quad \psi^p_N = \nabla_p\Phi(q_N,p_N),\label{eq:stateco6}\\
\Psi_{ki}^q & =  \psi^q_{k} + h \sum\limits_{j=1}^s \bar{a}_{ij}^q \phi^q(Q_{kj},P_{kj},\Psi_{kj}^q,\Psi_{kj}^p),\label{eq:stateco7}\\
\Psi_{ki}^p & = \psi^p_{k} + h \sum\limits_{j=1}^s \bar{a}_{ij}^p \phi^p(Q_{kj},P_{kj},\Psi_{kj}^q,\Psi_{kj}^p),\label{eq:stateco8}\\
\bar{a}_{ij}^q & = \frac{b_i b_j - b_j a^q_{ji}}{b_i},\quad \bar{a}_{ij}^p = \frac{b_i b_j - b_j a^p_{ji}}{b_i}\label{eq:stateco9}.
\end{align}
\end{subequations}
Assuming $h$ small enough and $x_k = (q_k,p_k)$ near $x^*(t_k) = (q^*(t_k),p^*(t_k))$ the intermediate states $Q_{ki}$ and $P_{ki}$ and costates $\Psi^q_{ki}$ and $\Psi^p_{ki}$ are uniquely determined since $u(x,\psi)$ depends Lipschitz continuously on $x$ near $x^*(t)$ and $\psi$ near $\psi^*(t)$ for any $t\in [0,T]$. 
This follows from smoothness and the implicit function theorem as for standard Runge-Kutta discretization as stated in \cite{Hager00}  (see for example \cite{Fla76}). 
Thus, there exists a locally unique solution $(Q_{ki},P_{ki},\Psi_{ki}^q, \Psi_{ki}^p),\, 1\le i \le s$ of \eqref{eq:stateco5}--\eqref{eq:stateco8}.

The scheme \eqref{eq:stateco} can be viewed as a discretization of the two-point boundary-value problem \eqref{eq:contbound}. To ensure a desired order of approximation for this two-point boundary-value problem, \cite{Hager00} derives order conditions on the coefficients of the Runge-Kutta scheme via Taylor expansions. In our case, however, we more easily obtain the following result concerning the order of consistency.

\begin{theorem}[Order of consistency]\label{th:orderapp}
If the symplectic partitioned Runge-Kutta discretization of the state system is of order $\kappa$ and $b_i>0$ for each $i$, then the scheme for the adjoint system is again a symplectic partitioned Runge-Kutta scheme of the same order (in particular we obtain the same schemes for $(q,p)$ and $(\psi^p,\psi^q)$).
\end{theorem}
\begin{proof} Starting with a symplectic partitioned Runge-Kutta discretization for the state system (\ref{eq:conv_ocp}), the discrete necessary optimality conditions are given as the adjoint system  (\ref{eq:stateco5}--\ref{eq:stateco9}). This system is again a partitioned Runge-Kutta scheme for the adjoint equations with coefficents $b_i, \bar{a}_{ij}^q, \bar{a}_{ij}^p$ . Substituting the symplecticity condition $a_{ij}^p = \frac{b_i b_j - b_j a^q_{ji}}{b_i}$ into equation (\ref{eq:stateco9}), the coefficients are determined as
\begin{align*}
&\bar{a}_{ij}^q = a_{ij}^p,\\
&\bar{a}_{ij}^p = a_{ij}^q.
 \end{align*}
 Since the coefficients of the Runge-Kutta scheme of the adjoint system for $(\psi^p,\psi^q)$ are the same as the coefficients of the Runge-Kutta scheme of the state system for $(q,p)$ defined in (\ref{eq:stateco1}--\ref{eq:stateco4}), the adjoint scheme is of same order.
\end{proof}

With Theorem \ref{th:ordercalculation} and \ref{equ_L_RK} one obtains the following:

\begin{lemma}\label{lemma_order}
Let a Lagrangian optimal control problem with regular Lagrangian $L$ and Lagrangian control force $f_{L}$ be given. If the discrete Lagrangian $L_d$ \eqref{eq:DMOC_RK_L} and the discrete control force $f_{d}^\pm$ \eqref{eq:DMOC_RK_f} with $b_i>0$ for each $i$ are both of order $\kappa$, then the corresponding adjoint scheme is also of order $\kappa$.
\end{lemma}

\subsection{Convergence}

The purpose of this section is to establish the convergence of  solutions of the discrete optimal control problem \eqref{eq:gDLOCP} to a solution of the Lagrangian optimal control problem \eqref{eq:LOCP} as the step size $h$ goes to $0$. 
%
%
As in Section \ref{adjoints}, we restrict ourselves to the case without path and control constraints and without final point constraint. Our convergence statement is a direct application of the one in \cite{DoHa00} (cf. also \cite{Hager00} and \cite{OB08}). 
In the previous section, we derived the adjoint system and determined the order of consistency. A standard strategy to show convergence of a discrete scheme is to prove consistency and stability. This is abstractly stated within the following result:
\begin{theorem}\cite{Hager01}\label{abstract_result}
Let $\mathcal{X}$ be a Banach space and let $\mathcal{Y}$ be a linear normed space with the norm in both spaces denoted by $\|\cdot\|$. Let $\mathcal{F}: \mathcal{X} \mapsto 2^\mathcal{Y}$ be a set-valued map, let $\mathcal{L}: \mathcal{X} \mapsto \mathcal{Y}$ be a bounded, linear operator, and let $\mathcal{T}: \mathcal{X} \mapsto \mathcal{Y}$ with $\mathcal{T}$ continuously Frech$\acute{\text{e}}$t differentiable in $B_r(w^*)$ for some $w^* \in \mathcal{X}$ and $r>0$. Suppose that the following conditions hold for some $\delta \in \mathcal{Y}$ and scalars $\epsilon, \lambda$ and $\sigma>0$:
\begin{itemize}
\item[{\rm (P1)}] $\mathcal{T}(w^*) + \delta \in \mathcal{F}(w^*)$.
\item[{\rm (P2)}] $\|\nabla \mathcal{T}(w) - \mathcal{L} \| \le \epsilon$ for all $w\in B_r(w^*)$.
\item[{\rm (P3)}] The map $(\mathcal{F} - \mathcal{L})^{-1}$ is single-valued and Lipschitz continuous in $B_\sigma(\pi), \pi = (\mathcal{T} - \mathcal{L}) (w^*)$, with Lipschitz constant $\lambda$.
\end{itemize}
If $\epsilon \lambda < 1$, $\epsilon r \le \sigma$, and $\| \delta\| \le (1-\lambda\epsilon) r/ \lambda$, then there exists a unique $w\in B_r(w^*)$ such that $\mathcal{T}(w) \in \mathcal{F}(w)$. Moreover, we have the estimate 
\begin{equation}\label{conv_erstimate}
\| w- w^*\| \le \frac{\lambda}{1-\lambda\epsilon} \|\delta\|.
\end{equation}
\end{theorem}
Consistency corresponds to assumption (P1) and the bounds on the norm of $\delta$, stability corresponds to assumption (P3) and the bound on the Lipschitz constant $\lambda$ for the linearization, and convergence is stated in (\ref{conv_erstimate}).

The following convergence result is formulated in terms of the averaged modulus of smoothness of the optimal control. If $J\subset\mathbb{R}$ is an interval and $v: J \rightarrow \mathbb{R}^n$, let $\omega(v,J;t,h)$ denote the modulus of continuity:
\begin{equation}
\omega(v,J;t,h) = \sup \{ | v(s_1)-v(s_2) |\, :\, s_1,s_2 \in [t-h/2, t+h/2] \cap J \}.
\end{equation}
The averaged modulus of smoothness $\tau$ of $v$ over $[0,T]$ is the integral of the modulus of continuity:
\[
\tau(v;h) = \int _0^T \omega (v,[0,T];t,h)\, \text{d} t.
\]
It is shown in \cite{SePo88} that $\lim_{h\to 0} \tau(v;h)=0$ if and only if $v$ is Riemann integrable, and $\tau(v;h) \le ch$ with constant $c$ if $v$ has bounded variation. 
\begin{theorem}[Convergence]\label{main_theorem}
Let $(x^*,u^*)$ be a solution of the Lagrangian optimal control problem \eqref{eq:LOCP1}--\eqref{eq:LOCP3}. If smoothness and coercivity hold, the discrete Lagrangian $L_d$ \eqref{eq:DMOC_RK_L} and the discrete control force $f_{d}^\pm$ \eqref{eq:DMOC_RK_f}  are of order $\kappa$ with $b_i>0$ for each $i$, then for all sufficiently small $h$ there exists a strict local minimizer $(x^h,u^h)$ of the discrete Lagrangian optimal control problem \eqref{eq:gDLOCP} and an associated adjoint variable $\psi^h$ satisfying the first-order necessary optimality conditions
such that
\begin{equation}\label{eq:conv_estimate}
\max\limits_{0\le k\le N} \left| x_k^h - x^*(t_k)\right| + \left| \psi_k^h - \psi^*(t_k)\right| + \left| u(x_k^h,\psi_k^h) - u^*(t_k)\right| \le ch^{\kappa-1} \left( h+ \tau \left( \frac{d^{\kappa-1}}{dt^{\kappa-1}} u^*;h \right) \right),  
\end{equation}
where $u(x_k^h,\psi_k^h)$ is a local minimizer of the Hamiltonian \eqref{eq:DMOC_Hamiltonian} corresponding to $x=x_k$ and $\psi = \psi_k$. 
\end{theorem}
 
\begin{proof} By Theorems \ref{th:ordercalculation} and \ref{equ_L_RK} we know that a discrete Langrangian and discrete forces both of order $\kappa$ lead to a symplectic partitioned Runge-Kutta discretization of order $\kappa$ for the state system. Because of smoothness and coercivity we can build up a discrete adjoint scheme with eliminated control that approximates the continuous adjoint scheme with order $\kappa$ (see Lemma \ref{lemma_order}). This leads in Hager's terminology to a {\em Runge-Kutta scheme of order $\kappa$ for optimal control}, and therefore Hager's convergence result for standard Runge-Kutta schemes in \cite{Hager00}, Theorem 2.1, is directly applicable. \end{proof}

\paragraph{Remark.}
Note that the estimate for the error in the discrete control in \eqref{eq:conv_estimate} is expressed in terms of $u(x_k^h,\psi_k^h)$ not $u_k$. 
This is due to the fact, that we derive the estimate via the transformed adjoint system with removed control due to the control uniqueness property. In \cite{DoHa00} the estimate is proved in terms of $u_k$ for Runge-Kutta discretization of second order.

\paragraph{Remark.}
Theorem \ref{main_theorem} can be extended to optimal control problems with constraints on the control function $u(t)$ as it was done in \cite{Hager00} for Runge-Kutta discretizations of order $2$.

\section{Implementation and Applications}

\subsection{Implementation}\label{impl}

As a balance between accuracy and efficiency we employ the midpoint rule for approximating the relevant integrals for the example computations in the following section, that is we set
\begin{eqnarray*}
C_d(q_k,q_{k+1},u_k)  & =& hC\left(\frac{q_{k+1}+q_k}{2},\frac{q_{k+1}-q_k}{h},u_{k+1/2}\right),\\
L_{d}(q_{k},q_{k+1}) &= &h L\left(\frac{q_{k+1}+q_k}{2},\frac{q_{k+1}-q_k}{h}\right),\\
\int_{kh}^{(k+1)h} f_{L}(q(t),\dot{q}(t),u(t))\cdot \delta q(t)\, \text{d} t &\approx& h f_{L}\left(\frac{q_{k+1}+q_k}{2},\frac{q_{k+1}-q_k}{h},u_{k+1/2}  \right)\cdot\frac{\delta q_{k+1}+\delta q_{k}}{2}\\
\end{eqnarray*}
$k=0,\dots,N-1$. Here, $f_{k}^-=f_{k}^+=\frac{h}{2} f_{L}\left(\frac{q_{k+1}+q_k}{2},\frac{q_{k+1}-q_k}{h},u_{k+1/2}  \right) $ are used as the left and right discrete forces with $q_k = q(t_k)$ and $u_{k+1/2}=u\left(\frac{t_k+t_{k+1}}{2}\right)$.

\paragraph{SQP Method.}
We solve the resulting finite dimensional constrained optimization problem by a standard SQP method as implemented for example in the routine  {\tt fmincon} of MATLAB. For more complex problems we use the routine {\tt nag\_opt\_nlp\_sparse} of the NAG library\footnote{\tt www.nag.com}. Since SQP is a local method, different initial guesses can lead to different solutions. This has been observed in almost all our example computations.
 
\paragraph{Automatic Differentiation.}

To compute an optimal solution of the optimization problem the SQP method makes use of the first and second derivatives of the constraints and the objective function. In the case where no derivatives are provided, the  algorithms that are used approximate those by finite differences. This approximation is time-consuming and round-off errors in the discretization process occur leading to worse convergence behavior of the algorithm. To avoid these drawbacks we make use of the concept of {\em Automatic Differentiation (AD)} (see \cite{Grie, Rall, Weng}), a method to numerically evaluate the derivative of a function specified by a computer program.  AD exploits the fact that any computer program that implements a vector function $y = F(x)$ (generally) can be decomposed into a sequence of elementary assignments, any one of which may be trivially differentiated by a simple table lookup. These elemental partial derivatives, evaluated at a particular argument, are combined in accordance with the chain rule from differential calculus to yield information on the derivative of $F$ (such as gradients, tangents, and the Jacobian matrix). This process yields (to numerical accuracy) exact derivatives. Because symbolic transformations occur only at the most basic level, AD avoids the computational problems inherent to complex symbolic computation. In particular, for some of our optimal control problems, we implemented the package ADOL-C ({\bf A}utomatic {\bf D}ifferentiation by {\bf O}ver{\bf L}oading in {\bf C}++, \cite{ADOL}) that has been written primarily for the evaluation of gradient vectors (rows or columns of Jacobians). It turns out that the optimization process performs in a faster and more robust way when providing the derivatives by automatic differentiation rather than via finite differences.

\subsection{Applications}\label{appl}

In this section, we numerically compare DMOC to a collocation method by means of two problems: a low thrust orbital transfer and the optimal control of a two-link manipulator. For this comparison, we apply a collocation method of order 2 to two different models: the Hamiltonian system with coordinates $(q,p)$ as well as the system formulated on tangent space with coordinates $(q,\dot{q})=(q,v)$.

\subsubsection{Low Thrust Orbital Transfer}\label{orb_trans}

In this application we investigate the problem of optimally transferring a satellite with a continuously acting propulsion system from one circular orbit around the Earth to another one.

\paragraph{Model.} Consider a satellite with mass $m$ which moves in the gravitational field of the Earth (mass $M$). The satellite is to be transferred from one circular orbit to one in the same plane with a larger radius, while the number of revolutions around the Earth during the transfer process is fixed.  In $2$d-polar coordinates  $q=(r,\varphi)$, the Lagrangian of the system has the form
\[
L(q,\dot{q}) = \frac{1}{2}m (\dot{r}^2+r^2\dot{\varphi}^2) + \gamma\frac{Mm}{r},
\]
where $\gamma$ denotes the gravitational constant. Assume that the propulsion system continuously exhibits a force $u$ in the direction of motion of the satellite, such that the corresponding Lagrangian control force is given by
\[
f_{L} = \left(\begin{array}{c}0\\ r\,u\end{array}\right).\]

\paragraph{Boundary Conditions.} Assume further that the satellite initially moves on a circular orbit of radius $r^0$.  Let  $(r(0),\varphi(0))=(r^0,0)$ be its position at $t=0$, then its initial velocity is given by $\dot{r}(0) = 0$ and $\dot{\varphi}(0) = \sqrt{\gamma M/(r^0)^3}$.  Using its thruster, the satellite is required to reach the point $(r^T,0)$ at time $T = d \sqrt{\frac{4\pi^2}{8\gamma M}(r^0+r^T)^3}$ and, without any further control input, to continue to move on the circle with radius $r^T$.  Here, $d$ is a prescribed number of revolutions around the Earth. Thus, the boundary values at $t=T$ are given by $(r(T),\varphi(T))=(r^T,0)$ and $(\dot{r}(T),\dot\varphi(T)) = (0,\sqrt{\gamma M/(r^T)^3})$.

\paragraph{Objective Functional.} During this transfer, our goal is to minimize the control effort, correspondingly the objective functional is given by
\[
J(q,\dot{q},u)=\int_0^T u(t)^2 \,\text{d} t.
\]
\paragraph{Results.} We compute the transfer from an orbit of radius $r^0=30$ km to one of radius $r^T=330$ km around the Earth.    First, we investigate how well the balance between the change in angular momentum and the amount of the control force is preserved. Due to the invariance of the Lagrangian under the rotation $\varphi$, the angular momentum of the satellite is preserved in the absence of external forces (as stated in Noether's theorem). However, in the presence of control forces, equation \eqref{eq:noether_force} gives a relation between the forces and the evolution of the angular momentum.

\begin{figure}[htbp]
\begin{center}
\includegraphics[width=0.48\textwidth]{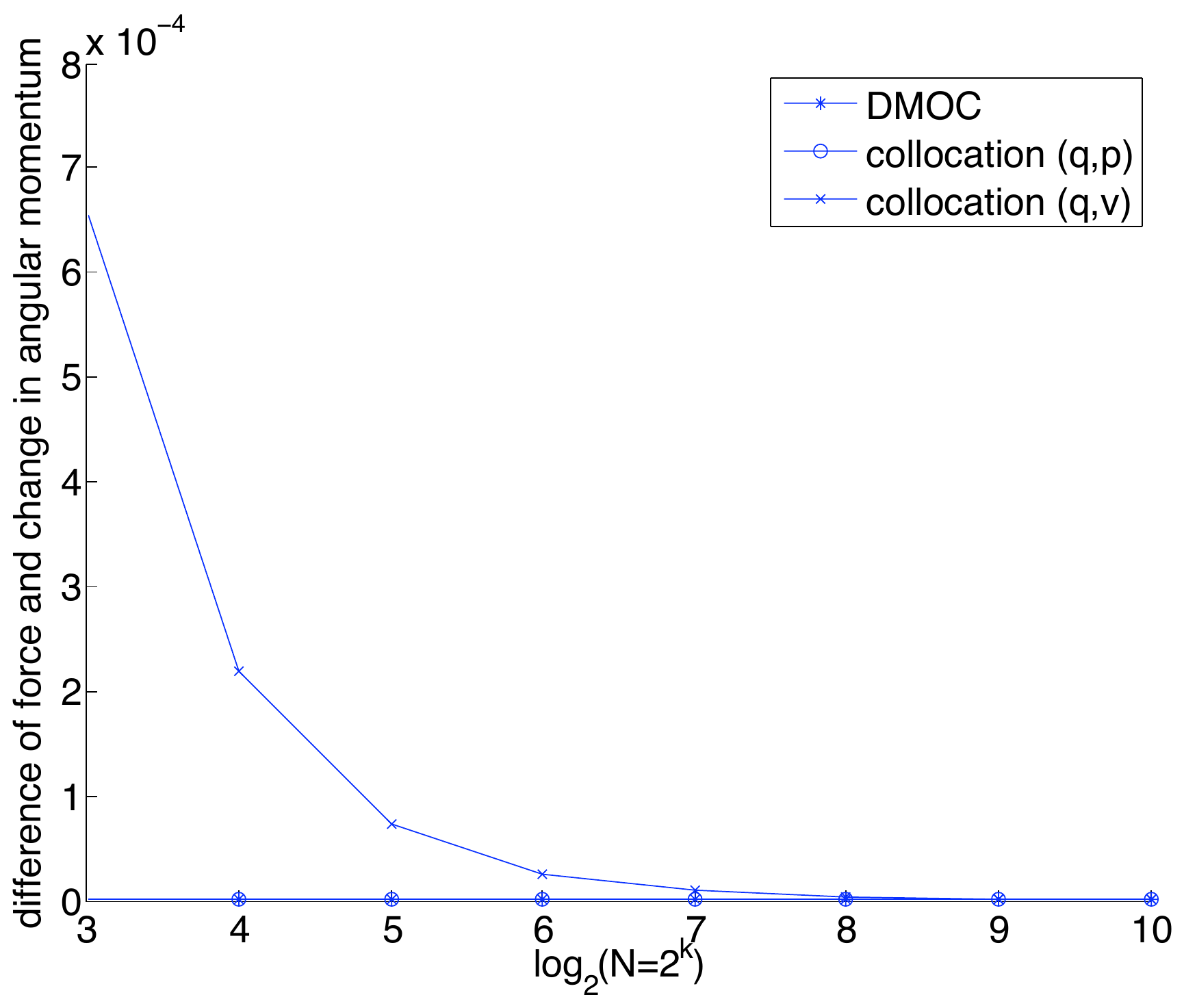}
\caption[Low thrust orbital transfer: difference of force and change in angular momentum]{Comparison of the accuracy of the computed open loop control for DMOC and a collocation approach:  difference of force and change in angular momentum in dependence on the number of discretization points. 
\label{fig:kepler_J}}
\end{center}
\end{figure}
In Figure~\ref{fig:kepler_J}, we compare the amount of the acting force with the change in angular momentum in each time interval. For the solution resulting from DMOC and the collocation approach applied to the Hamiltonian system, the change in angular momentum exactly equals the sum of the applied control forces (to numerical accuracy). The collocation method of second order corresponds to a discretization via the implicit midpoint rule. Thus, the optimization problem resulting from DMOC is equivalent to that obtained by applying collocation to the Hamiltonian formulation of the system as shown in Theorem~\ref{th:equi} and Example~\ref{ex:RK}a). Obviously, we obtain equal solutions by applying both methods. These results are consistent with the well-known conservation properties of variational integrators, that provide discretizations that preserve the continuous properties as momentum maps in the discrete setting in a natural way. On the other hand, the collocation method applied to the tangent space system described in velocities fails to capture the change in angular momentum accurately because the discrete tangent space formulation destroys the discrete Hamiltonian structure and the resulting (also unforced) scheme is not momentum-preserving anymore.

\begin{figure}[htbp]
\begin{center}
\includegraphics[width=0.48\textwidth]{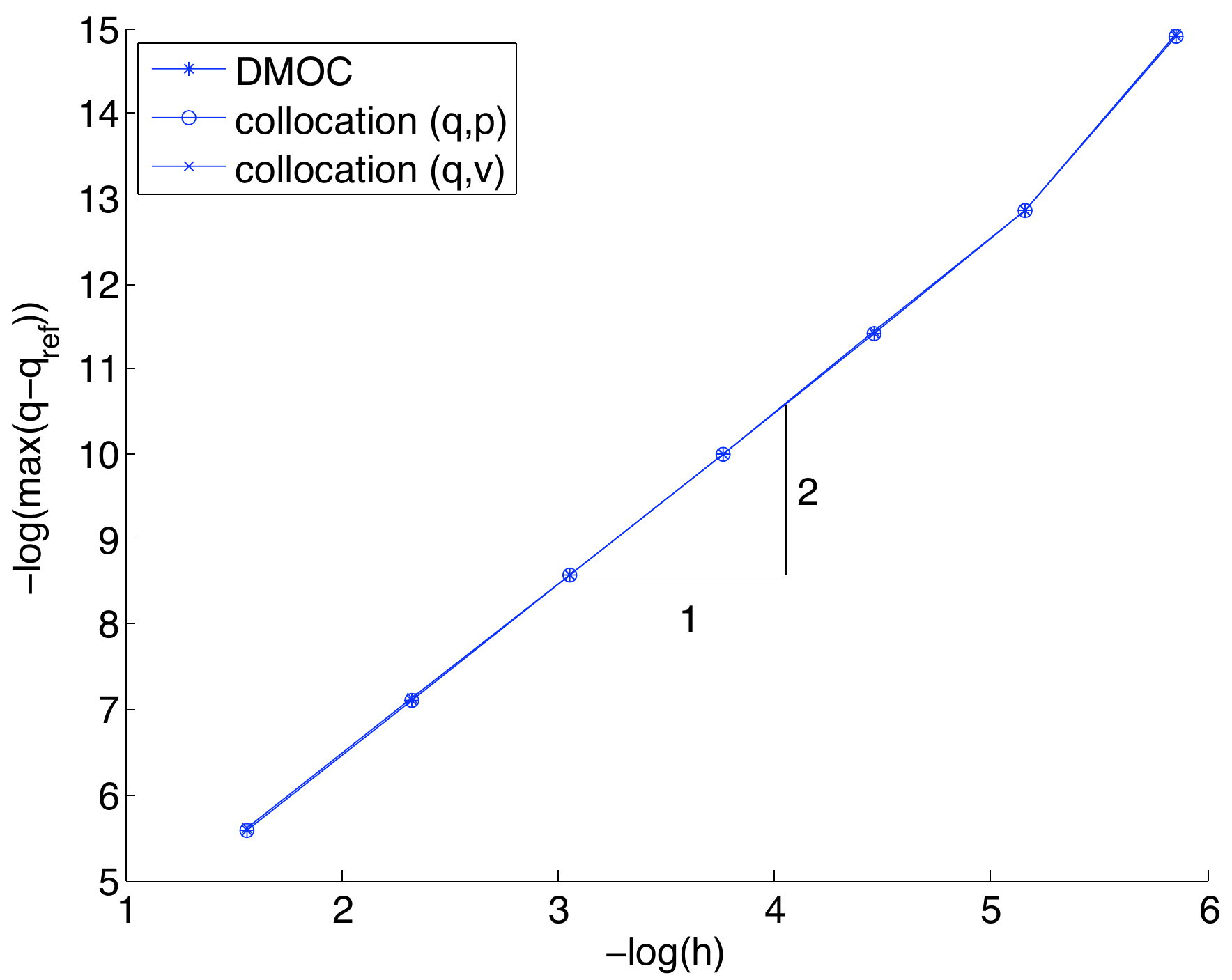}
\includegraphics[width=0.48\textwidth]{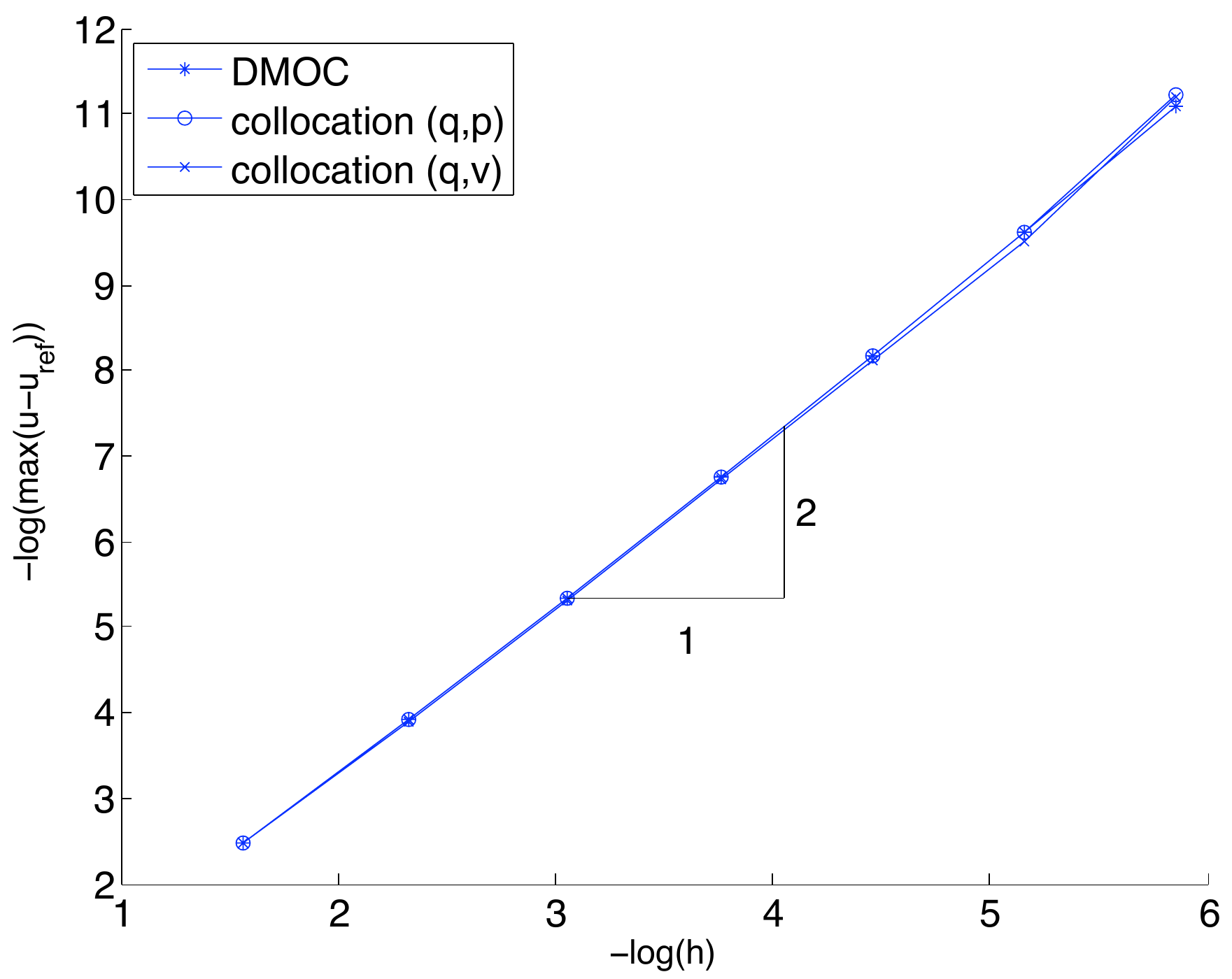}

\caption[Approximation error vs.\ step size]{Approximation error for the states (left) and the controls (right) in dependence of the step size \label{fig:kepler_error}}
\end{center}
\end{figure}
In Figure~\ref{fig:kepler_error} we show the convergence rates for all three methods. A reference trajectory is computed with $N=1024$ discretizations points and time step $h=2.9\cdot 10^{-3}$. The error in the configuration and control parameter of the discrete solution with respect to the reference solution is computed as 
$ \max_{k=0,\dots,N} |q(t_k)-q_\mathrm{ref}(t_k)|$ and $\max_{k=0,\dots,N} |u(t_k)-u_\mathrm{ref}(t_k)|$, respectively, where $|\cdot|$ is the Euclidean norm.
For all three methods the convergence rate for the configuration and control trajectory is $\mathcal{O}(h^2)$, as predicted by Theorem~\ref{main_theorem}.

\begin{figure}[htbp]
\begin{center}
\includegraphics[width=0.48\textwidth]{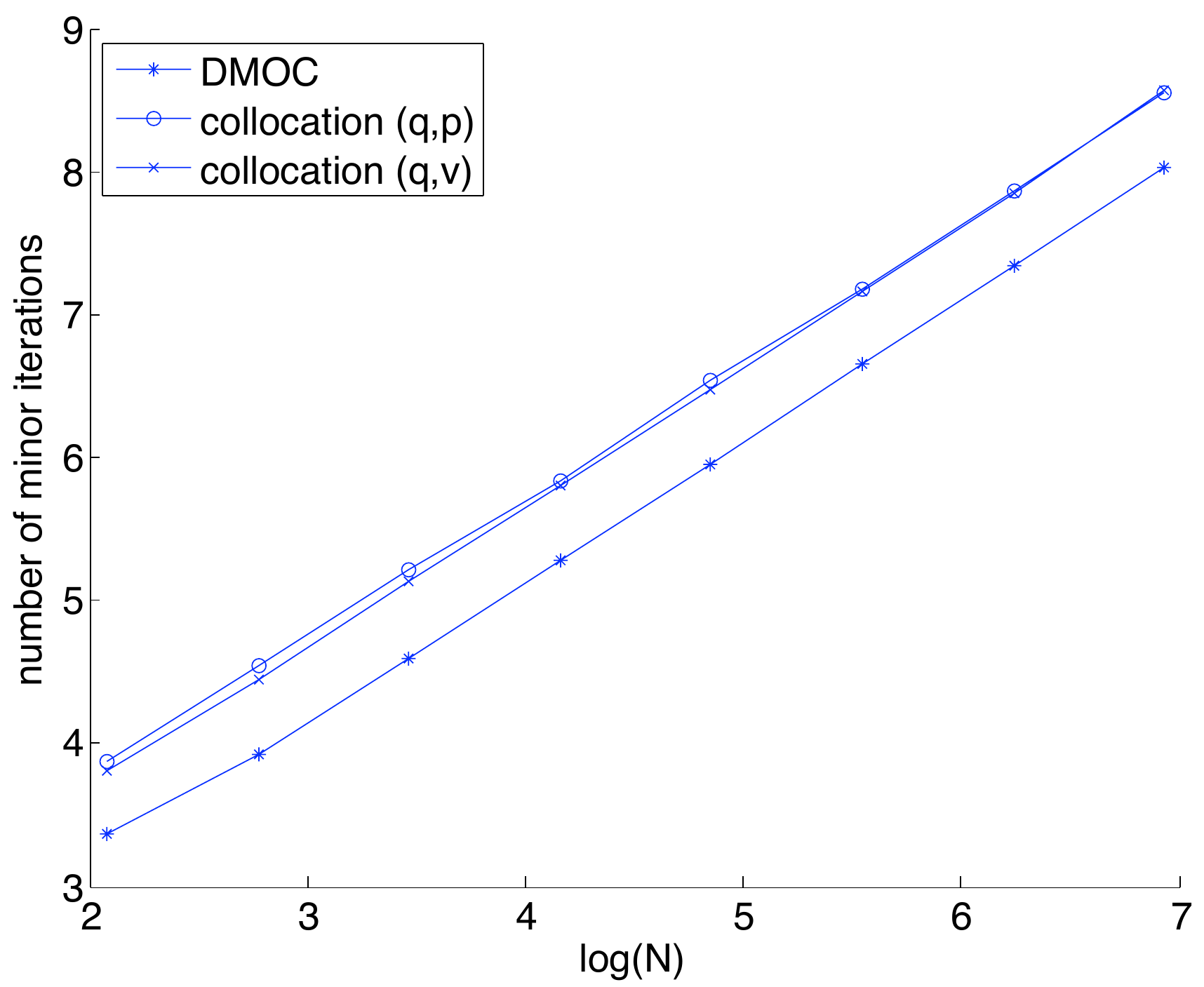}
\caption[Low thrust orbital transfer: Number of minor iterations]{Comparison of the number of iteration steps performed by the SQP solver for DMOC and a collocation approach in dependence on the number of discretization points.
\label{fig:kepler_iterations}}
\end{center}
\end{figure}
Still, DMOC is advantageous regarding the computational efficiency. Due to its formulation on the discrete configuration space, DMOC only uses $\frac{3}{5}\approx 0.6$ of the number of variables of the collocation approach (cf.~Remark~\ref{rem:savings}).  Figure~\ref{fig:kepler_iterations} shows the number of all iterations that the SQP solver performs in order to solve the quadratic subproblems (\emph{minor iterations}). We observe that for DMOC the SQP solver needs about $1.7$ times less iterations in comparison to collocation. This is in exact aggreement with the reduced number of variables.

\subsubsection{Two-link Manipulator}\label{imp_twolink}

As a second numerical example we consider the optimal control of a two-link manipulator. Again, we compare DMOC to  a collocation method of the same order. 

\begin{figure}[h!]
\begin{center}
\includegraphics[width=0.55\textwidth]{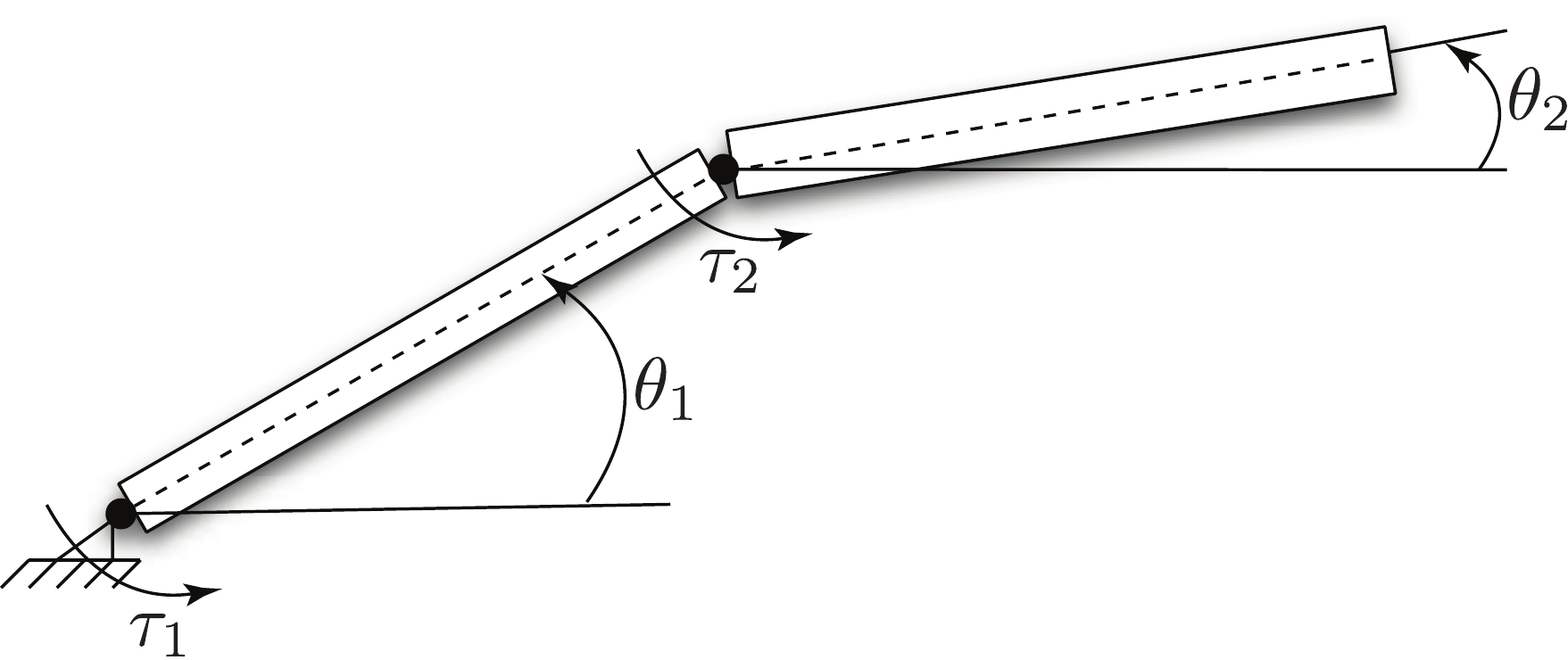}
\caption[Two-link manipulator: model]{Model of a two-link manipulator.}
\label{fig:twolink}
\end{center}
\end{figure}

\paragraph{Model.} The two-link manipulator (see Figure~\ref{fig:twolink}) consists of two coupled planar rigid bodies with mass $m_i$, length $l_i$ and moment of inertia $J_i$, $i=1,2$, respectively.
For $i\in {1,2}$, we let $\theta_i$ denote the orientation of the $i$th link measured counterclockwise from the positive horizontal axis. If we assume one end of the first link to be fixed in an inertial reference frame, the configuration of the system is specified by $q=(\theta_1,\theta_2)$. 

The Lagrangian is given via the kinetic and potential energy
\begin{eqnarray*}
K(q,\dot{q}) = \frac{1}{8}(m_1+4m_2)l_1^2\dot{\theta}_1^2 + \frac{1}{8} m_2 l_2^2\dot{\theta}_2^2 + \frac{1}{2}m_2 l_1 l_2 \cos{(\theta_1-\theta_2)}\dot{\theta}_1\dot{\theta}_2 + \frac{1}{2}J_1\dot{\theta}_1^2 + \frac{1}{2}J_2\dot{\theta}_2^2
\end{eqnarray*}
and
\begin{eqnarray*}
V(q) = \frac{1}{2}m_1gl_1 \sin{\theta_1} + m_2 g l_1 \sin{\theta_1} + \frac{1}{2} m_2gl_2 \theta_2,
\end{eqnarray*}
with the gravitational acceleration $g$. 
Control torques $\tau_1,\tau_2$ are applied at the base of the first link and at the joint between the two links. This leads to the Lagrangian control force
\[
f_{L}(\tau_1,\tau_2) = \left(\begin{array}{c} \tau_1-\tau_2 \\ \tau_2 \end{array}\right).
\]

\paragraph{Boundary Conditions.} The two-link manipulator is to be steered from the stable equilibrium point $q^0=(-\frac{\pi}{2}, -\frac{\pi}{2})$ with zero angular velocity $\dot{q}^0=(0,0)$ to the unstable equilibrium point $q^T=(\frac{\pi}{2}, \frac{\pi}{2})$ with velocity $\dot{q}^T=(0,0)$.

\paragraph{Objective Functional.}
For the motion of the manipulator the control effort 
\[
J(\tau_1,\tau_2) = \int _0^T \frac{1}{2}(\tau_1^2(t) + \tau_2^2(t))\, \text{d} t
\]
is to be minimized, where we fix the final time $T=1$.

\begin{figure}[htbp]
\begin{center}
\begin{tabular}{cc}
\includegraphics[width=0.48\textwidth]{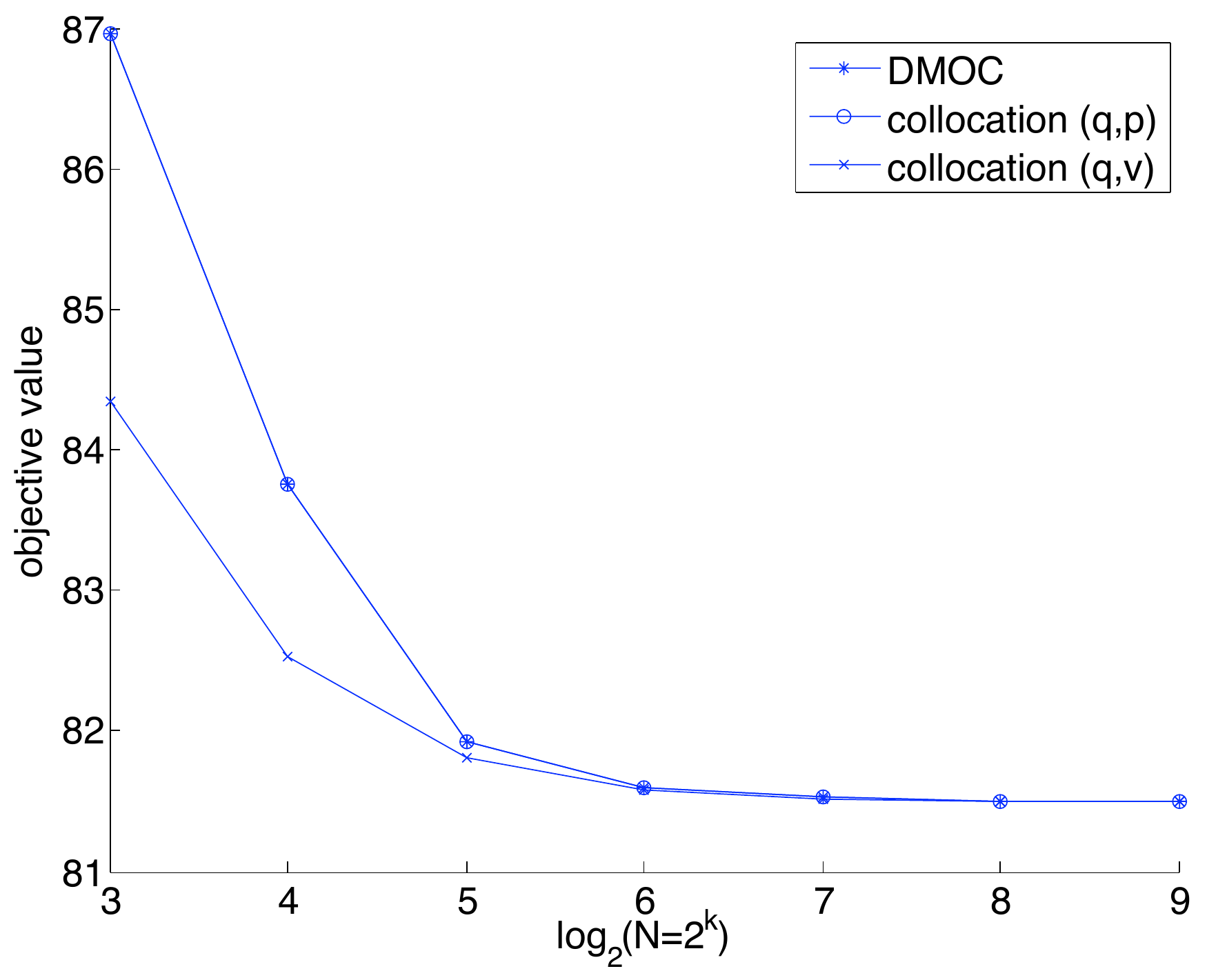}&
\includegraphics[width=0.48\textwidth]{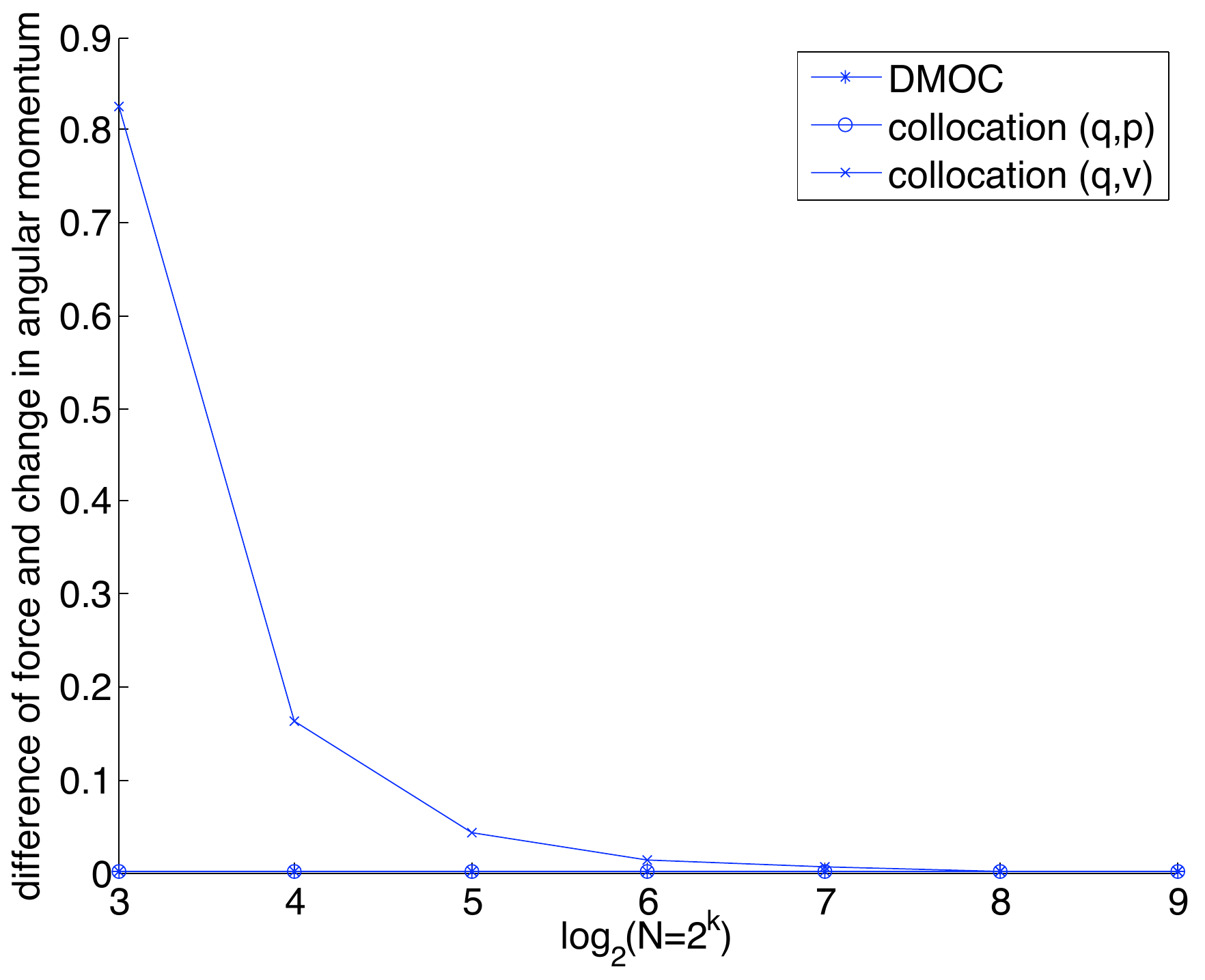} \\
{\footnotesize a)}&{\footnotesize b)}
\end{tabular}
\caption[Two-link manipulator: objective function value and momentum-force consistency]{Comparison of the accuracy of the computed open loop control for DMOC and a collocation approach: a) approximated cost b) difference of force and change in angular momentum in dependence on the number of discretization points.   \label{fig:twolink_obj}}
\end{center}
\end{figure}

\begin{figure}[htbp]
\begin{center}
\begin{tabular}{cc}
\includegraphics[width=0.48\textwidth]{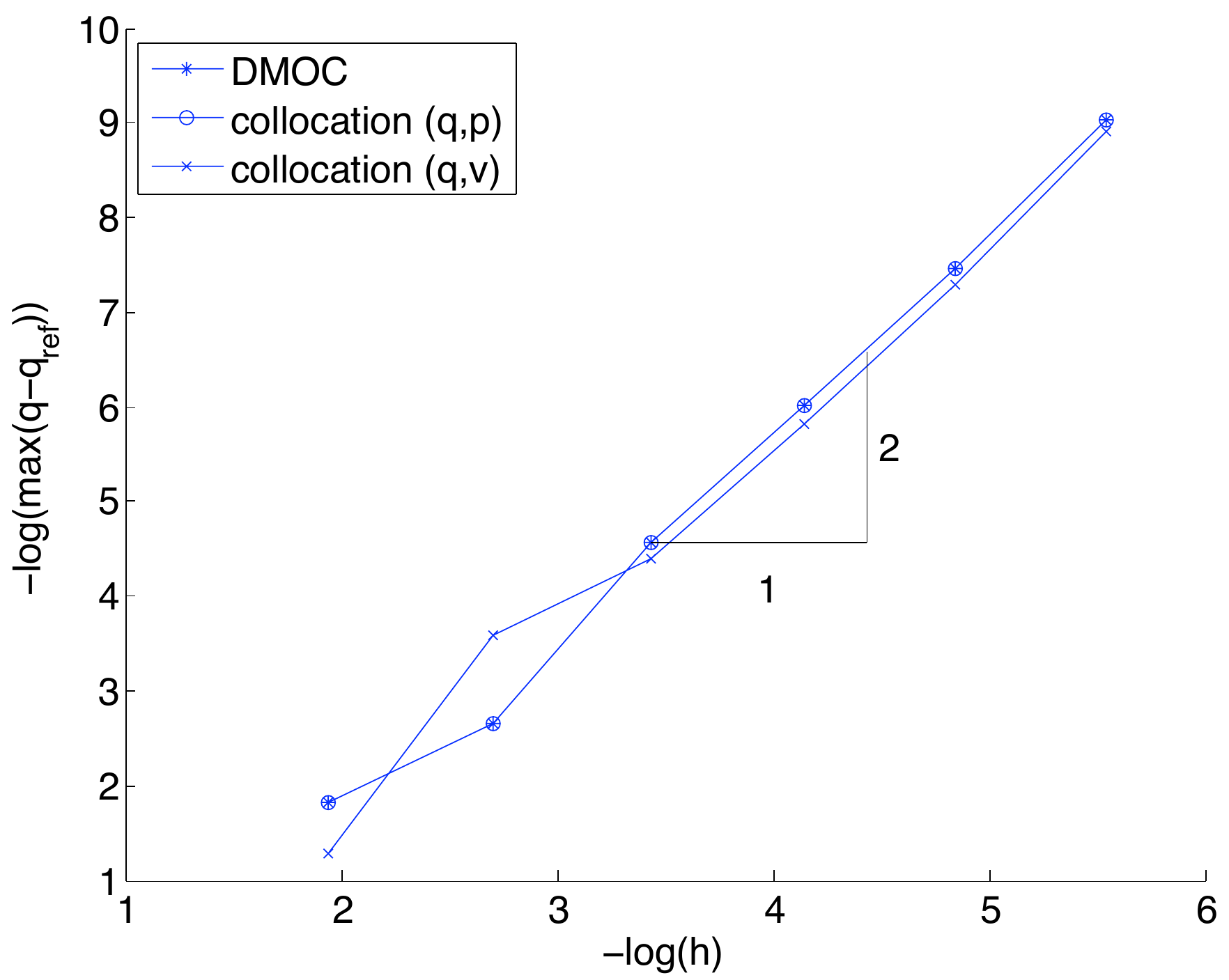}&
\includegraphics[width=0.48\textwidth]{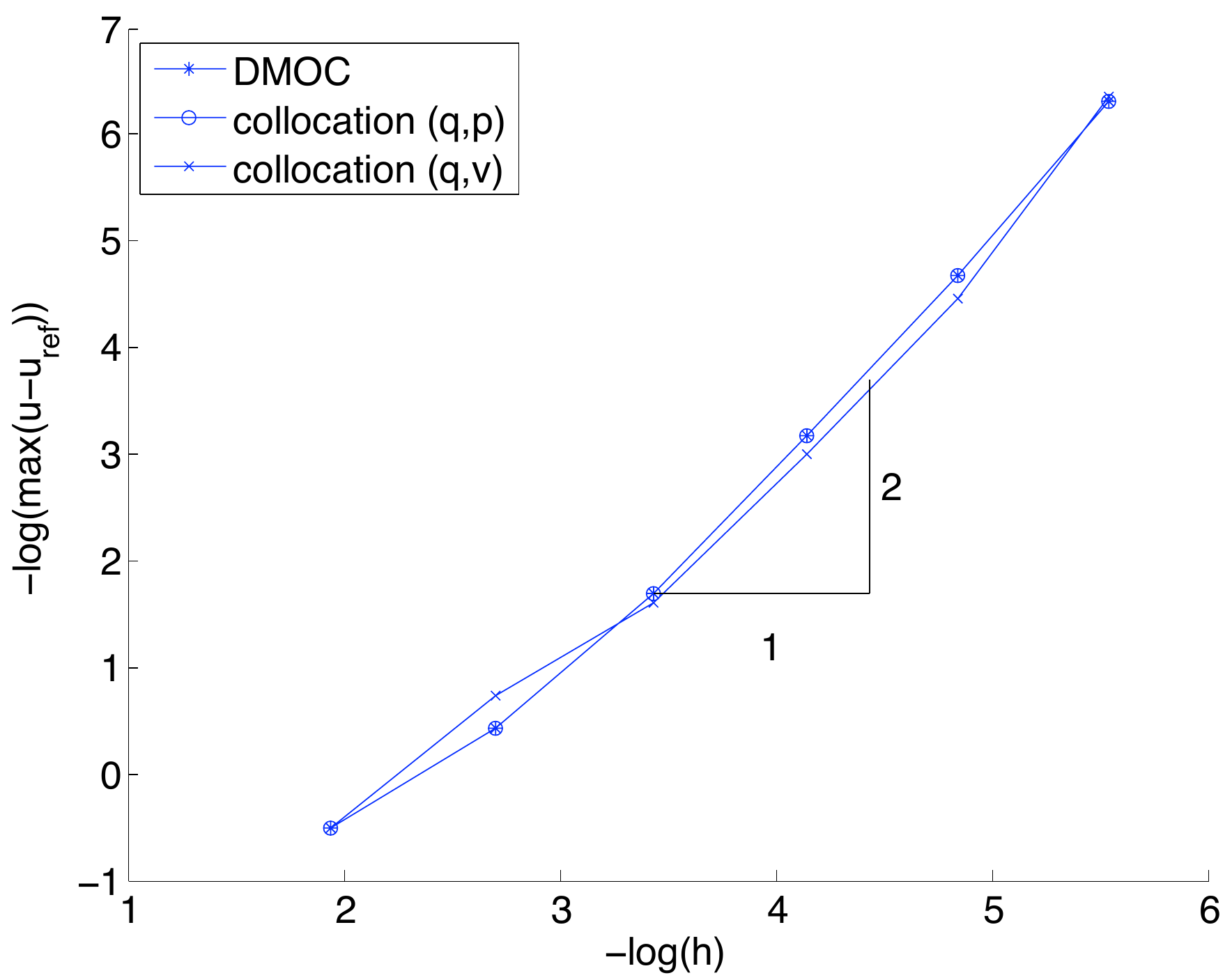} \\
{\footnotesize a)}&{\footnotesize b)}
\end{tabular}
\caption[Two-link manipulator: convergence rates]{Comparison of the convergence rates in dependence of the time step $h$ for DMOC and a collocation approach. a) Error of configuration trajectory. b) Error of control trajectory.   \label{fig:twolink_error}}
\end{center}
\end{figure}

\paragraph{Results.} In Figure~\ref{fig:twolink_obj} we show a) the resulting cost and b) the difference of the amount of force (including the control and the gravitational force) and the change in angular momentum in dependence on the number of discretization points for all three methods.  As expected, we obtain (numerically) identical solutions for DMOC and the equivalent collocation method for the Hamiltonian formulation.  The midpoint rule applied to the tangent space formulation performs equally well with respect to the objective value evolution. However, as in the previous example it does not reflect the momentum-force consistency as good as the other methods as shown in Figure~\ref{fig:twolink_obj}.   
In Figure~\ref{fig:twolink_error} the convergence rates are depicted. Here, a reference trajectory has been computed with $N=512$ discretizations points and time step $h=1.9\cdot 10^{-3}$. For all three methods the convergence rate for the configuration and control trajectory is $\mathcal{O}(h^2)$, as expected for a scheme of second order accuracy.

In Figure~\ref{fig:kepler_iterations} the number of minor iterations that the SQP solver performs to solve the quadratic subproblems is shown. Similar to what we have seen in the orbital transfer, we observe that the SQP solver needs about $1.5$ times less iterations when using DMOC rather than a collocation approach. This factor reflects the fact, that DMOC uses only $\frac{2}{3}\approx 0.67$ of the number of variables in comparision to the collocation approach.

\begin{figure}[htbp]
\begin{center}
\includegraphics[width=0.48\textwidth]{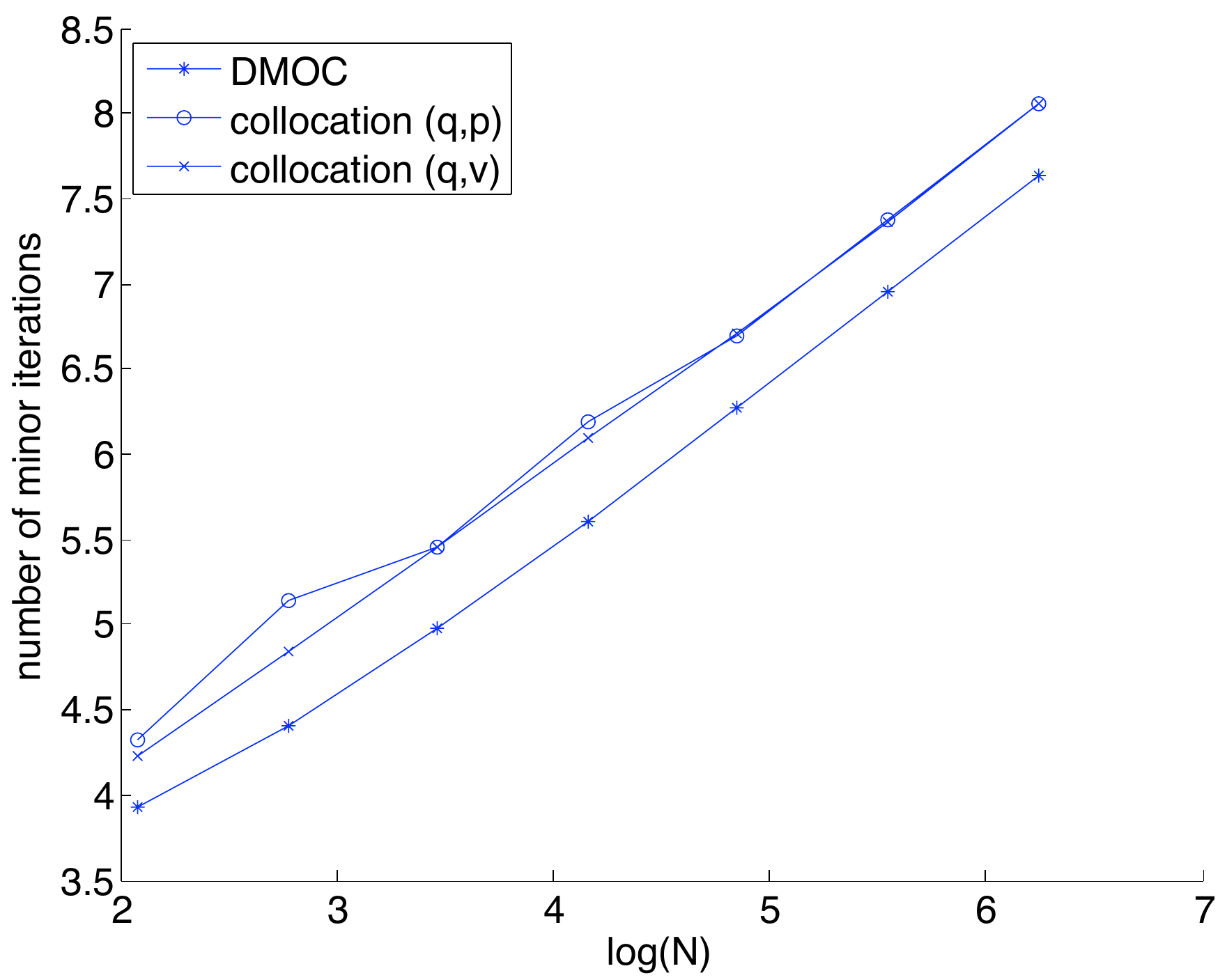}
\caption[Low thrust orbital transfer: Number of minor iterations]{Comparison of the number of iteration steps performed by the SQP solver for DMOC and a collocation approach in dependence on the number of discretization points.  \label{fig:twolink_iterations}}
\end{center}
\end{figure}

\section{Conclusions and Future Directions}

In this paper we developed a fully discrete formulation of the optimal control of mechanical systems. Based on discrete variational principles, the discrete optimal control problem yields a strategy, denoted by DMOC, for the efficient solution of these kinds of problems. 
The benefits of using DMOC are twofold: On the one hand, in the presence of a symmetry group in the continuous dynamical system, it is known that the momentum map changes in a definite way. The use of discrete mechanics allows one to find an exact counterpart to this on the discrete level. In this paper, this behavior was shown numerically for specific examples. On the other hand, due to the fact that DMOC is implemented on  the configuration level rather than on the configuration-momentum or configuration-velocity level, one gets significant computational savings for the number of iterations steps the SQP solver needs to solve the optimization problem. Again, this was demonstrated numerically in examples.

Although the developed method works very successfully in many examples considered in this paper and others (see \cite{JMaOb,JOb,JMOb,DMOCC,OB08}), there are still some open challenges that will be investigated in future work. 

\paragraph{Global Minima}
One challenge is the issue of local versus global minima. The use of an SQP solver for solving the resulting optimization problem restricts one to local minima   that are dependent on the specific initial guess used. One approach to overcome this issue is the use of DMOC primitives (\cite{KOB08,FDF05}). Here, the authors create a roadmap of feasible trajectories given by a graph where the edges represent small pre-computed DMOC segments referred to as DMOC primitives. A feasible solution of the optimal control problem corresponds to a specific concatenation of the DMOC primitives respecting the dynamics. The global optimal control can now be approximated by determining the optimal path in the graph with a given cost function  with a global search using techniques from dynamic programming. Initial demonstrations of this idea show that it can be carried out in real time. This, together with the fact that {\it real} dynamics is included via DMOC primitives are the significant advantages of this approach.

\paragraph{Adaptive time-stepping}

For many applications, an adaptive time-stepping strategy is essential. For example, for problems 
in space mission design, the planned trajectories require a finer time-stepping nearby planets due 
to the strong influence of gravity, while for a transfer in nearly free space only few discretization 
points are necessary to accurately reflect the dynamics of the system.
Here, different strategies such as error control based on the discretization grid under consideration and variational approaches could be investigated. For the variational approach a constraint is included to the Lagrange-d'Alembert principle that ensures time step control directly at the level of the discrete action (\cite{KMLTMD}). According to different adaption schemes different constraints can be included such as adaption to acceleration or the strenght of control forces.

\paragraph{Miscellaneous}
The framework could be extended to the optimal control of mechanical systems with stochastic influence or contact problems, respectively making use of the theory of stochastic (\cite{BRO07}) and nonsmooth variational integrators (\cite{FMOW03}), respectively. Due to the variational formulation of DMOC, an extension towards the optimal control of partial differential equations for the treatment of fluid and continuum mechanics might be interesting as well (see \cite{MS99,MPS98} for basic extensions of variational integrators to the context of PDEs).

\begin{small}
\bibliographystyle{marsden}
\bibliography{DMOCRefs}
\clearpage
\end{small}

\end{document}